\documentclass[a4paper,fleqn,sort&compress]{cas-sc}

\ExplSyntaxOn
\cs_set:Npn \__first_footerline: {}
\ExplSyntaxOff

\usepackage[numbers]{natbib}
\usepackage{stmaryrd}
\usepackage{amsmath}
\usepackage{amsfonts}
\usepackage[american]{babel}
\usepackage{subcaption}
\usepackage{hyperref}
\usepackage[utf8]{inputenc}
\usepackage{booktabs}
\DeclareUnicodeCharacter{2212}{−}
\usepackage{placeins}

\usepackage{mathtools}
\usepackage{amssymb}
\usepackage{amsthm}
\usepackage{enumitem}
\usepackage{bm}
\usepackage{bbm}
\usepackage{algorithm}
\usepackage[italicComments=false, rightComments=false]{algpseudocodex}
\usepackage{siunitx}
\usepackage{tikz,pgf}
\usetikzlibrary[patterns]
\usetikzlibrary{shapes,snakes}

\usepackage{xcolor}
\usepackage{ulem}
\usepackage{todonotes}

\usepackage[capitalise,nameinlink]{cleveref}
\newtheorem{rmk}{Remark}[section]
\newtheorem{defin}{Definition}[section]
\crefname{rmk}{Remark}{Remarks} 
\Crefname{rmk}{Remark}{Remarks} 
\crefname{defin}{Definition}{Definitions} 
\Crefname{defin}{Definition}{Definitions}

\DeclareMathOperator{\spn}{span}

\def\tsc#1{\csdef{#1}{\textsc{\lowercase{#1}}\xspace}}
\tsc{WGM}
\tsc{QE}
\tsc{EP}
\tsc{PMS}
\tsc{BEC}
\tsc{DE}

\newcommand{\modtibo}[1]{#1}
\newcommand{\corclem}[1]{#1}

\newcommand{\pablo}[2]{%
    {\color{red}#1}%
    \ifx&#2&%
    \else%
        {\color{brown}\sout{#2}}%
    \fi%
}

\begin{document}
\let\WriteBookmarks\relax
\def\floatpagepagefraction{1}
\def\textpagefraction{.001}

\shorttitle{Efficient Simulation of Nonlinear Lattice Structures}   

\shortauthors{C. Guillet \& al.}

\title[mode=title]{Efficient Fine-Scale Simulation of Nonlinear Hyperelastic Lattice Structures}

%

\author[1,2]{C. Guillet}

\author[3]{T. Hirschler}
\ead{thibaut.hirschler@utbm.fr}
\cormark[1]

\author[2]{P. Jolivet}

\author[4]{P. Antolin}

\author[5,6,7]{R. Bouclier}




\affiliation[1]{organization={Concace, Inria Center at the University of Bordeaux}, city={Talence}, postcode={33405}, country={France}}
\affiliation[2]{organization={Sorbonne Université, LIP6, CNRS},  
            city={Paris}, 
            postcode={75005}, 
            country={France} }
\affiliation[3]{
    organization={Universit\'e Marie et Louis Pasteur, UTBM, CNRS, ICB UMR 6303}, 
    city={Belfort}, 
    postcode={90010 },
    country={France}}
\affiliation[4]{organization={Ecole Polytechnique Federale de Lausanne, Institute of Mathematics},  
            city={Lausanne}, 
            postcode={1015}, 
            country={Switzerland} }
\affiliation[5]{organization={Univ. Toulouse, INSA-Toulouse, IMT-Albi, ISAE-SUPAERO, CNRS UMR 5312, ICA}, 
            city={Toulouse}, 
            postcode={31400}, 
            country={France} }
\affiliation[6]{organization={Univ. Toulouse, INSA-Toulouse, CNRS UMR 5219, IMT}, 
            city={Toulouse}, 
            postcode={31077}, 
            country={France} }
\affiliation[7]{organization={Institut Universitaire de France (IUF)},
            country={France} }

\cortext[cor1]{Corresponding author}

\begin{abstract}
With the growing maturity of additive manufacturing, the fabrication of architected or lattice-based metamaterials has become a reality for industrial applications. These materials combine lightweight design with tailored mechanical properties, most of which exhibit pronounced nonlinear, especially large-deformation, behaviors. The main numerical challenge therefore lies in performing nonlinear simulations of such lattice structures, which may contain thousands of geometrically intricate unit cells, while lacking sufficient scale separation for multiscale homogenization schemes to be applicable straightforwardly. In this work, we propose a dedicated solver for the full volumetric fine-scale simulation of nonlinear hyperelastic lattice structures that drastically reduces both memory and computational costs. The key idea is to exploit the intrinsic self-similarity of the cells through a reduced-order modeling strategy applied within a domain-decomposition framework. At each Newton iteration, a limited set of principal cells is identified through a dedicated, weakly intrusive, EIM-like approach, allowing all local tangent operators to be expressed as linear combinations of a few principal ones. This enables fast and memory-efficient operator assembly, and then feeds an efficient inexact FETI-DP–based preconditioner at the solution stage, resulting in a quasi matrix-free algorithm for the nonlinear analysis. Numerical experiments in two and three dimensions demonstrate significant computational gains, with runtime reductions from several hours to a few tens of minutes and memory savings by factors of about three, while maintaining full fine-scale accuracy. Notably, the proposed strategy enables the computation of problems involving thousands of cells (i.e., millions of degrees of freedom) within a few minutes on an off-the-shelf laptop.
\end{abstract}


\maketitle


\section{Introduction}

With the advent of additive manufacturing, which, as its name suggests, produces parts by adding rather than removing material, the fabrication of industrially relevant architected materials has become a reality~\cite{blakey21,benedetti21}. A particularly effective class of such materials is what can be referred to as lattice structures or lattice microstructures, featuring both carefully designed unit-cell geometries and optimized macroscopic shapes~\cite{antolin19optim,weeger22iga,zwar23shape,ntop,antolin2025design}. The periodic (or quasi periodic) architectures of these materials enable the combination of lightweight design with tailored mechanical properties, resulting not only in high stiffness-to-weight ratios but also in novel behaviors such as high stretchability, auxetic effects, and multifunctional responses. Owing to these unique features, lattice structures are often regarded as mechanical metamaterials, which are highly attractive for a wide variety of engineering applications, including thermal management, mechanical resistance, and energy absorption in aerospace, automotive, civil, and defense engineering, as well as in the biomedical field~\cite{vasiliev12,schaedler16,tancogne16,bastek23inverse,wu23additively}.

In this work, we focus on the numerical simulation of the mechanical response of lattice structures. In view of the final goal of applying our strategy to realistic lattice structures, we consider the large-deformation regime and, at this stage, adopt the framework of hyperelasticity. This modeling approach is particularly well suited to accurately capture the mechanical response of polymer lattice structures under representative loading conditions~\cite{fernandez21,nakarmi24,dos25size}. However, this task poses a major computational challenge, as industrially relevant lattice structures can now contain thousands of geometrically intricate unit cells. Consequently, the computational cost in terms of both memory and time becomes extremely demanding, or even intractable, when standard approaches such as the finite element method are used as black-box solvers. This explains why almost all volumetric fine-scale simulations of these structures reported in the literature are restricted to a limited number of unit cells. Only a few studies seem to push this limit further by relying on massively parallel high-performance computing architectures, but at the expense of substantial computational resources~\cite{korshunova21bending,shaikeea22,brown22,de24}. Therefore, there is a growing need to develop more efficient numerical strategies that reduce the computational burden at its source, enabling the accurate prediction of the structural response of lattice structures within reasonable computational budgets, an aspect that becomes even more critical in nonlinear regimes.

To address this issue, the most common trend is to act at the modeling level in order to circumvent the direct fine-scale simulation of such structures. In other words, the full volumetric discretization of the architectural geometry, i.e., at the scale of the struts or walls composing the unit cells, is avoided. Among these strategies, one can first mention multiscale homogenization methods, which were initially developed for general heterogeneous materials and have since been enhanced through data-driven techniques~\cite{gartner21,korshunova21bending,somnic22,masi22,glaesener23}. In such approaches, only fine-scale problems of the size of one or a few cells need to be solved. The cells are actually regarded as “materials”: the full fine-scale geometrical details are replaced by equivalent constitutive laws that reproduce the averaged mechanical response of the cells. These methods have also become the approach of choice for integration within optimization loops, for instance to identify unit cells with prescribed effective properties~\cite{bastek23inverse,thakolkaran25} or to optimally populate a macroscopic domain with preselected cell topologies~\cite{costa22,liu22}. However, it should be noted that such approaches rely on a clear scale separation between the unit-cell size and the macroscopic dimensions. This assumption is often questionable in the case of lattice structures~\cite{yoder18,dos25size}, since 3D printing technologies constrain the achievable cell size, and the macro-shape frequently exhibits a slender geometry, thereby leading to only a few cells along certain directions. Consequently, another class of approaches consists in modeling lattice structures as networks of beams~\cite{weeger22iga,passieux23,glaesener23} or shells~\cite{ma24}, which allows representing the full architecture of the lattice structure while reducing the computational complexity. Yet, there are inherent challenges associated with the use of such models for lattice structures, including the proper incorporation of complex nonlinear constitutive laws~\cite{weeger22mixed} and the accurate connection between different structural components (struts and/or walls)~\cite{portela18,cadart25}. Moreover, this type of modeling may not be suitable for all cell geometries, especially those involving thick members.

In this work, we propose to follow a rather different conceptual path. The idea is to exploit the unique characteristics of lattice structures, namely their periodic (or quasi-periodic) nature, not at the modeling level, but instead in the preconditioning of an iterative solver dedicated to computing the full volumetric fine-scale mechanical response of the structure. In other words, the goal is to solve the complete fine-scale problem while reducing the memory and computational costs by leveraging the (quasi) self-similarity of the cells throughout the lattice structure. Such an idea has recently begun to emerge in the literature, highlighting the potential of combining model reduction techniques, used to extract the common dominant mechanical modes of the cells, with domain-decomposition-based algorithms for preconditioning the problem~\cite{Hirschler2022,Hirschler:2023aa,GUILLET2025114136,rubio25}. Such frameworks have demonstrated the ability to achieve highly efficient matrix assembly and solver strategies, enabling, in particular, the prediction of the linear elastic response of lattice structures with millions of degrees of freedom on standard workstations~\cite{Hirschler:2023aa}, and up to billions of degrees of freedom on massively parallel architectures~\cite{GUILLET2025114136}. Building upon these previous contributions, the present work aims to extend such approaches to the large-deformation regime, considering at this stage a hyperelastic constitutive model. This extension is intended to generalize and consolidate these methods for the simulation of lattice structures under more representative operating conditions. 

The key challenge therefore lies in enabling a computationally efficient extraction of the common dominant mechanical modes of the cells as they undergo large deformations during the nonlinear simulation. This is achieved through the development of a dedicated reduced-order modeling (ROM) strategy, which can be related to the Empirical Interpolation Method (EIM) or its discrete counterpart, the DEIM~\cite{maday2013generalized,quarteroni15,hesthaven16}. The ROM is applied at each Newton iteration both to accelerate the assembly process and to precondition the inexact FETI-DP algorithm of~\cite{Hirschler:2023aa}, thereby leading to fast and memory-efficient solutions of the tangent systems. The resulting method is generic, in the sense that it can be applied to any type of discretization (standard finite elements or others), and remains weakly intrusive, as it operates at the discrete level, provided that one can assemble cell-wise operators using both reduced and full integration rules. As for the discretization scheme adopted in this work, we rely, for enhanced performance, on a Computer-Aided Design (CAD) paradigm based on spline composition~\cite{Elber_2017,elber23review} to generate the lattice geometries. Such a geometrical representation offers great flexibility for the design of lattice structures, since both the local (cell-level) and global (macro-scale) geometries are naturally parameterized by splines. It is also fully consistent with the IsoGeometric Analysis (IGA) framework~\cite{hirschler19embedded,MASSARWI2018148,Hirschler2022,GUILLET2025114136}, thereby facilitating the use of higher-order and smoother approximation spaces for the mechanical solution. Consequently, regarding the application of the developed approach, the present work also contributes to extending the state of the art in efficient IGA (see, e.g.,~\cite{marussig23fast,montardini2023ieti,montardini2025low,GUILLET2025114136} for recent contributions) toward nonlinear analysis.

The paper is organized as follows. In \cref{sec:1}, we first introduce all the necessary ingredients to set up the context of this study, namely the geometric modeling aspects and the nonlinear large-deformation hyperelastic framework. Then, in \cref{sec:discretization}, we present the developed methodology for fast operator assembly and efficient nonlinear analysis of lattice structures. Although the method is described within the context of IGA discretization, remarks are provided to highlight that it can be applied to more standard finite element schemes as well. Finally, a series of numerical experiments involving different cell patterns and macroscopic shapes, in both two and three dimensions, are presented in \cref{sec:res_num} to assess the efficiency of the proposed solver in terms of memory usage and computational cost reduction, before concluding remarks are drawn in \cref{sec:conclu}.

\section{Lattice structure modeling and nonlinear hyperelastic problem}
\label{sec:1}

This section introduces the context of the study and the corresponding notations. First, we briefly describe the mathematical model used to define the lattice structures, which is based on the composition of B-spline or NURBS functions. Then, we present the formulation of the nonlinear hyperelasticity problem.

\begin{figure}
    \centering
    \includegraphics[width=0.85\textwidth]{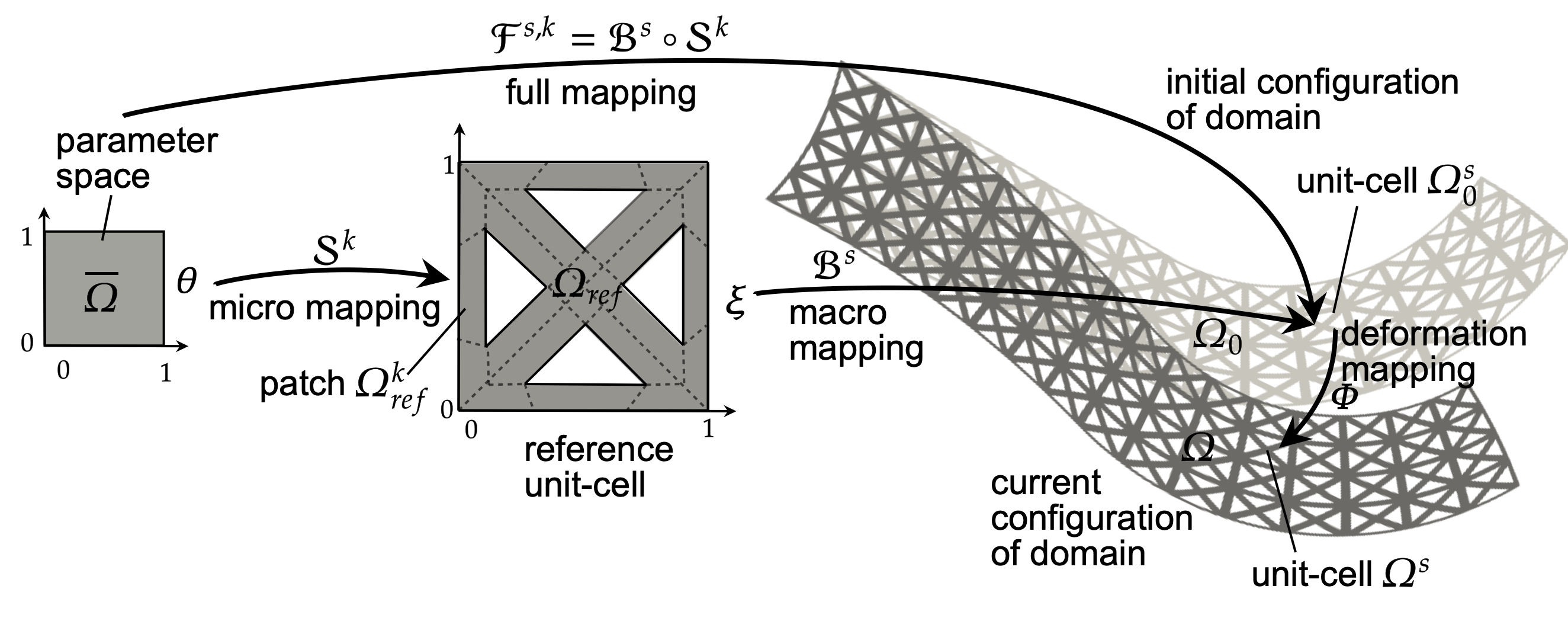}
	\caption{Brake pedal lattice structure, inspired from~\cite{ntopo} and defined by spline composition of a reference microstructure unit-cell and a macro-geometry mapping \label{fig:2}}
\end{figure}

\subsection{Modeling of lattices structures}
\label{sec:1.2}

In this work, lattice structures are defined and modeled at the fine-scale (or microstructural) level, i.e., where each individual cell is explicitly meshed. More precisely, this is achieved using a functional composition approach based on splines (either B-splines or NURBS), as introduced in \cite{Elber_2017,elber23review}. 

\corclem{
First, we briefly introduce the spline model, and refer the reader for more details to the books~\cite{cottrell09iso, piegl2012nurbs, bouclier2022iga} and the references therein. Given a set of control points with coordinates $\bm{P}_{\bm{i}}\in\mathbb{R}^d$, we defined a $d$-dimensional B-spline entity by
\begin{align}
\label{eq:1}
\mathcal{S}(\bm{\xi}) = \sum_{i=1}^{n} \phi_{i}^{\bm{p}}(\bm{\xi}) \bm{P}_{i}, ~~ \mathrm{with} ~~ \phi_{\bm{i}}^{\bm{p}}(\bm{\xi}) = \prod_{j=1}^{d}\phi_{i_j}^{p_j}(\bm{\xi}_j), ~~ \mathrm{for}~~ \bm{\xi} \in \bar{\Omega}:=[0,1]^d. \end{align}
In the above equation, $\phi_{i_j}^{p_j}$ denote the $n_j$ univariate $p_j$-degree B-spline basis functions, for $j\in \llbracket1,d\rrbracket$, and $n=\prod_{j=1}^{d}n_j$ denotes the total number of basis functions. To obtain multivariate NURBS functions and entities, it is necessary to adopt a rational formulation. In practice, this is achieved by associating a weight with each control point and defining the NURBS functions as the ratio of the weighted B-spline functions to the sum of all weighted B-spline functions. In the remainder of the paper, we will refer to \cref{eq:1} as a general single-patch spline mapping, representing either a B-spline or a NURBS mapping. Note that for complex geometries, it is necessary to use multi-patch spline models, i.e., models constructed by combining several patches. 
}

The spline composition approach consists of two components, illustrated in \cref{fig:2}(initial configuration): a macro-representation of the lattice structure, in which heterogeneities are not explicitly represented, and a reference microstructure that defines the pattern to be tiled within the macro-geometry to generate the cells. The final heterogeneous structure is obtained by embedding the reference microstructure into the macro-model through functional composition.

The first component in constructing the lattice structures is a reference unit cell, built by assembling multiple spline patches. 
\corclem{We introduce a reference unit cell, denoted $\Omega_{\mathrm{ref}}$ and represented in \cref{fig:2}, which is constructed from reference patches using a multi-patch spline model following the general spline mappings in \cref{eq:1} as 
\begin{align}
\Omega_{\mathrm{ref}} := \bigcup_{k=1}^{N_p}\Omega_{\mathrm{ref}}^{k}, \quad \mathrm{where}~~ \Omega_{\mathrm{ref}}^{k}:=\mathcal{S}^{k}(\bar{\Omega}), \quad \mathrm{and}~~ \mathcal{S}^{k}(\bm{\theta}) = \sum_{i=1}^{n_p} \phi_{i}^{k}(\bm{\theta}) \bm{P}_{i}^{k}  ~~\text{for}~ \bm{\theta} \in \bar{\Omega}.
\label{eq:micro}
\end{align}
These reference patches are assumed to be fully matching; that is, the intersections $\Omega_{\mathrm{ref}}^{k} \cap \Omega_{\mathrm{ref}}^{l}$, for $k,l \in\llbracket1,N_p\rrbracket$ such that $k\neq l$, are either empty, common vertices, common edges, or common faces (in 3D), and the discretizations on the common boundaries of these patches match perfectly. Here and throughout the following, we assume that all patches share the same degree. Accordingly, we omit the superscript $\bm{p}$ in the basis functions indicating the degree. Finally, $n_p$ denotes the number of basis functions associated with a given patch.
}

%
The second component is a macroscopic model representing the global shape of the lattice structure, prescribed here by a second spline model. This macro-mapping naturally induces a partition of the domain into $N_s\in\mathbb{N}^*$ macro-elements, each of which will contain one cell in the final geometry. Following the procedure in \cite{Hirschler2022}, we rewrite this macroscopic model using Bézier extraction so that the parameter spaces associated with each macro-element become identical. Consequently, the only difference between macro-elements lies in their underlying geometric mappings. 
\corclem{
The Bézier macro-mappings associated with the macro-elements are given by
\begin{align*}
\mathcal{B}^{s}(\bm{\xi}) = \sum_{i=1}^{\tilde{n}} \psi_{i}(\bm{\xi}) \bm{Q}_{i}^{s}, \quad \mathrm{for}~~\bm{\xi} \in \bar{\Omega}_M, \quad \forall s\in\llbracket1,N_s\rrbracket,
\end{align*}
where $\bar{\Omega}_M=[0,1]^d$ is the parametric space of the macro-model, $\psi_{i}$ denotes the Bernstein basis functions and $\bm{Q}_{i}^{s}$ stands for the associated control points.
}

\corclem{
The lattice structure domain (in its initial configuration), denoted by $\Omega_0$, is defined by the composition of the micro- and macro-mappings and consists of the union of $N_s$ fully-matching unit cells, denoted by $\Omega_0^{s}$, each of which being decomposed into $N_p$ patches. Specifically, it reads as
\begin{align*}
\Omega_0 = \bigcup_{s=1}^{N_s} \bigcup_{k=1}^{N_p}\Omega_0^{s,k}, \quad  \mathrm{where}~~\Omega_0^{s,k}:=\left\{ {\bf X} \in \mathbb{R}^d \ | \ {\bf X} = \mathcal{F}^{s,k}(\bm{\theta}), \quad \forall \bm{\theta} \in \bar{\Omega} \right\}.
\end{align*}
For all cells $s\in \llbracket1,N_s\rrbracket$ and patches $k\in\llbracket 1,N_p\rrbracket$, the unit-cell patches $\Omega^{s,k}_0$ are generated with the composition of the reference micro-model and the Bézier macro-mappings:
\begin{align*}
\mathcal{F}^{s,k} \colon & \bar{\Omega} \to \Omega_0^{s,k} \nonumber \\
& \bm{\theta} \mapsto (\mathcal{B}^{s} \circ \mathcal{S}^{k})(\bm{\theta}).
\end{align*}
}
%
Once again, see \cref{fig:2}(initial configuration) for an illustrative representation of this spline composition model.

\subsection{Nonlinear hyperelastic equations}
\label{sub:hyperelastic}

In this work, we consider the framework of nonlinear hyperelasticity to model the mechanical response of lattice structures. Two sources of nonlinearity are addressed and combined in this modeling: material and geometrical nonlinearities. The former occurs when the constitutive law, which relates stress to strain, is inherently nonlinear. The latter arises from the nonlinear relationship between displacement and strain, which becomes significant when the geometry of the body changes considerably under loading. Accounting for both types of nonlinearity is essential for accurately capturing the mechanical response, notably that of polymer lattice structures when they undergo large deformations. In the following, we detail the hyperelastic formulation adopted in this work, as it is necessary for understanding the developments presented in \cref{sec:discretization}. The presentation is based in particular on~\cite{bonet1997nonlinear}.

\subsubsection{Finite deformations}  
We begin by introducing the description of motion and the notations required for the analysis of finite deformations. Let $\Omega_0$ denote the initial (undeformed) configuration of the domain, i.e., before any forces are applied, and let $\Omega$ denote the current (deformed) configuration, i.e., after the application of forces (see \cref{fig:2} for an illustration). All previously introduced notations are extended to both the initial and current configurations. In this work, we adopt a spatial, or Eulerian, description of the body’s behavior; that is, all relevant quantities are expressed with respect to the current, deformed configuration. It is worth noting that a material, or Lagrangian, description could be equivalently adopted, leading to the same physical predictions.

The spatial position of a material point $\mathbf{x}\in\Omega$ is related to its initial position $\mathbf{X}\in \Omega_0$ through the deformation mapping $\bm{\Phi}$:
\begin{align*}
\mathbf{x} = \bm{\Phi}(\mathbf{X}) := \mathbf{X} + \bm{u}(\mathbf{X}), \quad \text{for}~\mathbf{X} \in \Omega_0,
\end{align*}
where $\bm{u}(\mathbf{X})$ is the displacement field. It follows that the current domain is the image of the initial domain under the deformation mapping: $\Omega = \bm{\Phi}(\Omega_0)$. 

The deformation gradient tensor is defined as the Jacobian of the deformation mapping:
\begin{align}
\label{eq:def_grad}
\bm{F}(\mathbf{X}) := \frac{\partial \bm{\Phi}}{\partial \mathbf{X}}(\mathbf{X}) = \mathbf{I} + \frac{\partial \bm{u}}{\partial \mathbf{X}}(\mathbf{X}), \quad \text{for}~\mathbf{X} \in \Omega_0,
\end{align}
where $({\bf I})_{ij} = \delta_{ij}$ is the identity tensor, and $\delta_{ij}$ is the Kronecker delta defined by
\begin{align}
\label{eq:kron}
\delta_{ij} = \left\{
\begin{array}{ll}
1 & \mathrm{if}~~ i = j, \\
0 & \mathrm{otherwise}.
\end{array}
\right.
\end{align}
Note that, here and throughout the remainder of this paper, the dependence of $\Omega$ and $\bm{F}$ on the deformation mapping~$\bm{\Phi}$ (or equivalently on the displacement field $\bm{u}$) is omitted for notational simplicity.

\subsubsection{Weak form of the equilibrium problem}

The boundary of the domain is denoted by $\partial \Omega$. A portion of the boundary, denoted by $\partial \Omega_D$, is assumed to be clamped, i.e., subject to homogeneous Dirichlet boundary conditions. The remaining part, denoted by $\partial \Omega_N = \partial \Omega \setminus \partial \Omega_D$, is subjected to a prescribed traction $\bm{g}$, which may vanish on a portion of $\partial \Omega_N$. In addition, a body force $\bm{f}$, such as gravity, is considered.

The appropriate functional setting for the variational formulation of the equilibrium problem is the Sobolev space $\bm{H}^1_0(\Omega, \partial \Omega_D) := \{\bm{v} \in \bm{H}^1(\Omega) \ | \ \bm{v} = 0 \text{ on } \partial \Omega_D\}$, where $\bm{H}^1(\Omega) = [H^1(\Omega)]^d$ stands for the $d$-dimensional Sobolev space. The subscript $()_0$ indicates here enforcement of Dirichlet boundary conditions, and should not be confused with the notation for the initial configuration. According to the principle of virtual work, the weak form of the equilibrium problem for a deformable body reads:
\begin{align}
\label{eq:virt_eq}
\text{Find } ~ \bm{\Phi} ~ \text{ s.t. } ~ W(\bm{\Phi}, \bm{v}) = 0, \quad \forall \bm{v} \in \bm{H}^1_0(\Omega, \partial \Omega_D).
\end{align}
This problem consists in finding the equilibrium configuration of the body by determining the deformation mapping $\bm{\Phi}$, which describes the transformation from the initial to the deformed configuration. Finding $\bm{\Phi}$ is equivalent to computing the displacement field $\bm{u}$.

The virtual work functional can be decomposed as
\begin{align*}
W (\bm{\Phi}, \bm{v}) = W_{\mathrm{int}} (\bm{\Phi}, \bm{v}) - W_{\mathrm{ext}} (\bm{\Phi}, \bm{v}),
\end{align*}
where the internal and external contributions are given by
\begin{align*}
W_{\mathrm{int}} (\bm{\Phi}, \bm{v}) = \int_{\Omega} \bm{\sigma}(\bm{\Phi}) : \delta \bm{d}(\bm{v}) \, d{\bf x}, \quad \mathrm{and}~~
W_{\mathrm{ext}} (\bm{\Phi}, \bm{v}) = \int_{\Omega} \bm{f}(\bm{\Phi}) \cdot \bm{v} \, d{\bf x} + \int_{\partial \Omega_N} \bm{g}(\bm{\Phi}) \cdot \bm{v} \, d{\bf s}.
\end{align*}
Here, $\bm{\sigma}$ denotes the Cauchy stress tensor, which depends on the nonlinear material behavior under consideration, and $\delta \bm{d}$ is the symmetric part of the virtual rate of deformation, defined by
\begin{align*}
\delta \bm{d}(\bm{v}) = \frac{1}{2} \left( \nabla \bm{v} + \nabla \bm{v}^T \right), \quad \text{for }~ \bm{v} \in \bm{H}^1_0(\Omega, \partial \Omega_D).
\end{align*}
In this work, we assume that the external virtual work is independent of the deformation, which is the standard case when considering gravity for $\bm{f}$ and non-follower forces for $\bm{g}$. Under this assumption, the external work can be equivalently computed in the initial configuration as
\begin{equation}
    W_{\mathrm{ext}}(\bm{v}) 
    = 
    \int_{\Omega_0} \bm{f}_0 \cdot \bm{v} \, d\mathbf{x} + \int_{\partial \Omega_{N,0}} \bm{g}_0 \cdot \bm{v} \, d{\bf s},
    \label{eq:externalvirtwork}
\end{equation}
where $\bm{f}_0$ and $\bm{g}_0$ are the applied loading in the initial configuration.

\subsubsection{Constitutive equations}
\label{sec:const_eq}
The virtual work~\cref{eq:virt_eq} is expressed in terms of the internal stresses induced by the strain of the material. The relationship between strain and the resulting stress is governed by the constitutive equations. In this work, the constitutive behavior is defined within the framework of a compressible hyperelastic neo-Hookean material, where stresses are derived from a strain energy density function.

The Cauchy stress tensor for a neo-Hookean material is given by
\begin{align}
\label{eq:sig_eq}
\bm{\sigma}(\bm{\Phi}) = \frac{\mu}{\det(\bm{F})} (\bm{b}(\bm{\Phi}) - \bm{I}) + \frac{\lambda}{\det(\bm{F})} \ln(\det(\bm{F})) \bm{I},
\end{align}
where $\bm{b}(\bm{\Phi})= \bm{F} \bm{F}^T$ is the left Cauchy–Green deformation tensor. 

We introduce the linearized strain tensor $\bm{\varepsilon}$ as
\begin{align*}
\bm{\varepsilon}(\bm{u}) = \frac{1}{2}\left(\bm{\nabla} \bm{u} + \bm{\nabla} \bm{u}^{T}\right).
\end{align*}
The stress–strain relationship for a neo-Hookean material can be expressed as
\begin{align*}
    D_{\bm{u}} \bm{\sigma}(\bm{\Phi}) =\bm{\mathcal{C}}(\bm{\Phi}) : \bm{\varepsilon}(\bm{u}),
\end{align*}
where $D_{\bm{u}} \bm{\sigma}(\bm{\Phi})$ denotes the directional derivative, measuring the variation of $\bm{\sigma}(\bm{\Phi})$ when $\bm{\Phi}$ is perturbed in the direction of $\bm{u}$ (this operator will be introduced in the next section).  
In the above expression, $\bm{\mathcal{C}}$ is the fourth-order hyperelasticity tensor, defined in the current configuration as
\begin{align}
\label{eq:const_eq}
\bm{\mathcal{C}}(\bm{\Phi}) 
= 
{\mathcal{C}}_{ijkl}(\bm{\Phi})\, \mathbf{e}_i \otimes \mathbf{e}_j \otimes \mathbf{e}_k \otimes \mathbf{e}_l,
\end{align}
where the components are given by
\begin{align*}
{\mathcal{C}}_{ijkl}(\bm{\Phi}) = \frac{\lambda}{\det(\bm{F})} \delta_{ij} \delta_{kl} + \frac{\mu - \lambda \ln(\det(\bm{F}))}{\det(\bm{F})} \left(\delta_{ik} \delta_{jl} + \delta_{il} \delta_{jk}\right).
\end{align*}
Here, $\delta$ is the Kronecker delta defined in \cref{eq:kron}, and $\mathbf{e}_i \in \mathbb{R}^d$ denotes the $i$-th canonical basis vector. The parameters $\lambda$ and $\mu$ are the Lamé coefficients, expressed with Poisson's coefficient $\nu$ and Young modulus $E$ as
\begin{align*}
\lambda = \frac{E\nu}{(1+\nu)(1-2\nu)} \quad \mathrm{and}~~ \mu = \frac{E}{2(1+\nu)}.
\end{align*}
Note here that the Cauchy stress tensor $\bm{\sigma}$ and the hyperelasticity tensor $\bm{\mathcal{C}}$ depend on the deformation mapping $\bm{\Phi}$, or equivalently, on the displacement field $\bm{u}$. 

\begin{rmk}
The neo-Hookean material model is adopted here for the sake of illustration. 
However, the methodology developed in \cref{sec:discretization} can be readily applied to other hyperelastic constitutive models, such as the Mooney–Rivlin formulation.
\end{rmk}

\subsubsection{Linearized equilibrium equations}

Problem \cref{eq:virt_eq} is nonlinear with respect to both geometry and material behavior. Geometric nonlinearities arise because the domain depends on the deformation mapping, which itself depends on the solution of the equilibrium problem. Material nonlinearities originate from the hyperelastic constitutive law~\cref{eq:sig_eq}, which is inherently nonlinear.

To compute the equilibrium configuration of the body, i.e., to solve problem \cref{eq:virt_eq}, a Newton–Raphson iterative method is employed. This requires linearizing the associated weak form. To do so, we compute the directional derivative of the virtual work functional with respect to a displacement increment $\bm{u}$ applied to the deformation mapping $\bm{\Phi}$. This directional derivative, denoted by $D_{\bm{u}} W (\bm{\Phi}, \bm{v})$, measures the change in $W (\bm{\Phi}, \bm{v})$ when $\bm{\Phi}$ is perturbed in the direction of $\bm{u}$, and is formally defined as:
\begin{align}
D_{\bm{u}} W (\bm{\Phi}, \bm{v}) := \lim_{h \rightarrow 0} \frac{W (\bm{\Phi} + h\bm{u}, \bm{v}) - W (\bm{\Phi}, \bm{v})}{h}.
\label{eq:direc_deriv}
\end{align}
Since the external work is assumed to be independent of the deformation, its directional derivative vanishes. Consequently, the linearization of the virtual work functional reduces to the internal work component, whose directional derivative reads:
\begin{equation}
    \begin{split}
    D_{\bm{u}} W (\bm{\Phi}, \bm{v})&= D_{\bm{u}} W_\mathrm{int} (\bm{\Phi}, \bm{v}) \\
    &= \int_{\Omega} \delta {\bf d} (\bm{v}) : \bm{\mathcal{C}}(\bm{\Phi}) : \varepsilon(\bm{u}) d{\bf x} + \int_{\Omega} \bm{\sigma}(\bm{\Phi}) : (\bm{\nabla u}^T \bm{\nabla v})  d{\bf x}.
    \end{split}    
    \label{eq:dir_der}
\end{equation}
The two terms of \cref{eq:dir_der} are commonly referred to as the constitutive and initial stress components, respectively.

\newcommand{\supiter}[1]{^{(#1)}}
\newcommand{\currentiter}{\mathsf{i}}
\newcommand{\nextiter}{\mathsf{i+1}}
\newcommand{\supcurrentiter}{\supiter{\currentiter{}}}
\newcommand{\supnextiter}{\supiter{\mathsf{i+1}}}
\newcommand{\uNewton}[1]{\mathbf{u}\supiter{#1}}
\newcommand{\PhiNewton}[1]{\bm{\Phi}^{h\mathsf{(#1)}}}

\subsubsection{Newton–Raphson method}

\modtibo{%


Let $\bm{\Phi}\supiter{0}$ be an initial guess for the deformation mapping, with associated displacement field $\bm{u}\supiter{0}$ satisfying:
\begin{align*}
    \bm{\Phi}\supiter{0}(\mathbf{X}) = \mathbf{X} + \bm{u}\supiter{0}(\mathbf{X}), \quad \forall \mathbf{X} \in \Omega_0.
\end{align*}
Given an approximation $(\bm{u}\supcurrentiter{}, \bm{\Phi}\supcurrentiter{})$ at iteration $\currentiter$, the Newton–Raphson update is defined by
\begin{align*}
    \bm{u}\supnextiter{} 
    = 
    \bm{u}\supcurrentiter{} 
    - \left(D_{\bm{u}} W (\bm{\Phi}\supcurrentiter{}, \bm{v})\right)^{-1} W(\bm{\Phi}\supcurrentiter{},\bm{v}), 
    \quad \forall \bm{v}\in \bm{H}^1_0(\Omega, \partial \Omega_D).
\end{align*}
The iteration process stops when the Newton update satisfies the nonlinear problem~\cref{eq:virt_eq} \corclem{to a prescribed tolerance}.

In a discrete setting, a linear system, commonly referred to as the tangent system, is solved at each iteration to update the solution and obtain an improved approximation. The precise discretization of the problem using IGA will be detailed later on in \cref{sec:discretization}. For now, let us assume an arbitrary discretization with a finite-dimensional subspace $\bm{V}_h \subset \bm{H}^1_0(\Omega, \partial \Omega_D)$, let us denote by~$\uNewton{\currentiter}$ the vector of displacement Degrees Of Freedom (DOF) at iteration~$\currentiter{}$. To compute the Newton increment
, we generate a descent direction~$\Delta\uNewton{\currentiter}$ by solving the following linear system
\begin{equation}
    \mathbf{K}(\uNewton{\currentiter}) \Delta\uNewton{\currentiter}
    =
    -\mathbf{r}(\uNewton{\currentiter}),
    \label{eq:linsys}
\end{equation}
where the tangent stiffness matrix~$\mathbf{K}$ and the residual load vector~$\mathbf{r}$ depend on the current Newton iteration. This might lead to high computational cost as it requires reassembly of the matrix and solution of the corresponding linear system at each iteration. To facilitate the convergence of the Newton method, a line search procedure is introduced to determine an optimal step length~$\alpha\supcurrentiter{} \in \mathbb{R}$ and to update the current approximation as
\begin{equation*}
    \uNewton{\nextiter} 
    = 
    \uNewton{\currentiter}
    +
    \alpha\supcurrentiter{}
    \Delta\uNewton{\currentiter}.
\end{equation*}
We adopt the so-called backtracking line search method, also known as Armijo's rule, to determine a suitable value for $\alpha\supcurrentiter{}$~\cite{zienkiewicz05finite,Nocedal2006-vy}. Specifically, we seek the smallest integer $m \geq 0$ such that
\corclem{\begin{equation*}
    \lVert\mathbf{r}(
        \uNewton{\currentiter}
        +
        \beta^m
        \Delta\uNewton{\currentiter}
        )
    \rVert_2
    \leq 
    (1 - c \beta^m )
    \lVert\mathbf{r}(\uNewton{\currentiter})
    \rVert_2,
\end{equation*}}
where $\beta = \frac{1}{2}$ is a fixed reduction factor, and $c = 10^{-4}$ controls the sufficient decrease condition. Once~$m$ in hand, the step length is defined has~$\alpha\supcurrentiter{} = \beta^m$. The line search enhances the robustness of Newton iterations, particularly in highly nonlinear regimes or near bifurcation points, by providing an update that safely ensures a reduction of the residual. Finally, the method stops when the relative Euclidean norm of the residual vector meets the condition ${\|\mathbf{r}(\uNewton{\nextiter})\|_2}/{\| \mathbf{r}(\uNewton{0}) \|_2} \leq 10^{-6}$.
}

\section{Dedicated discretization, operator assembly, and solver for efficient analysis}
\label{sec:discretization}

Having established the general framework of our study, we now turn to the proposed methodology for the efficient numerical simulation of nonlinear hyperelastic lattice structures. The presentation is organized as follows. In \cref{sub:dis_IGA}, we introduce the discretization of the problem using IGA, detailing the construction of the discrete operators, with particular emphasis on the tangent matrix. In \cref{sub:fast_assembly}, we describe the ROM-based fast assembly strategy developed to reduce both the computational cost and the memory requirements. Finally, \cref{sub:efficient_solver} presents the implemented linear inexact domain decomposition solver, which complements the fast assembly strategy to ensure efficient solutions at each Newton iteration of the nonlinear algorithm.

\modtibo{%
For clarity and conciseness of notation throughout this section, we restrict ourselves to a single Newton iteration. Accordingly, the superscript denoting the iteration index is omitted. 
}

\subsection{IGA discretization}
\label{sub:dis_IGA}

Starting from a geometry represented by splines, IGA~\cite{hughes05iga,cottrell09iso,bouclier2022iga} relies on the isoparametric element concept, i.e., the same high-degree and smooth basis functions are employed both for the description of the computational domain and for the discretization of the unknown fields. This unified representation provides proper geometry modeling, by bringing closer CAD-based design and numerical simulation, and offers higher continuity across element interfaces compared to classical finite element methods.

In this work, we adopt the strategy proposed in~\cite{hirschler19embedded,MASSARWI2018148,Hirschler2022,GUILLET2025114136}, where the micro-model spline functions (see \cref{eq:micro}) are employed to define the discretization space. These micro-mapping splines are preferred to the spline basis obtained from the full composition mapping for two main reasons: (i) they prevent the unnecessarily high polynomial degrees induced by composition~\cite{van18crossing,Hirschler2022}, which inevitably increase computational cost, and (ii) they guarantee a uniform discretization of the displacement field across all unit cells, a crucial feature to enable efficient matrix assembly and solver strategies that exploit the repetitive character of lattice cells~\cite{Hirschler2022,Hirschler:2023aa,GUILLET2025114136}.

\subsubsection{Trial and test functions, and approximation subspace}

For the mechanical analysis, it is necessary to consider a functional space that is at least  $\mathcal{C}^0$ over the entire computational domain. As stated above, this space is derived from the micro-model spline functions $\left( \phi_{\corclem{a}}^k \right)$ defined in~\cref{eq:micro}. More precisely, we construct below the corresponding $\mathcal{C}^0$ basis functions by strongly gluing together matching basis functions across patch and cell interfaces. For illustration of the different mappings mentioned below, we refer the reader again to \cref{fig:2}.

\corclem{
We introduce the space of $\mathcal{C}^0$-continuous basis functions associated with the reference cell, defined by
\begin{align*}
\spn\big\{  \phi^{\mathrm{ref}}_{\corclem{a}} ~ | ~  {\corclem{a}}\in \llbracket 1,n_{\mathrm{ref}}\rrbracket \big\}=\big\{ 
        v \in \mathcal{C}^0(\Omega_{\mathrm{ref}})
        ~|~
        \forall{k}\in\llbracket1,N_p\rrbracket,~ v_{|\Omega_{\mathrm{ref}}^{k}} \in \spn\{  \phi_{\corclem{a}}^{k}\circ(\mathcal{S}^{k})^{-1}  ~ | ~ {\corclem{a}}\in \llbracket1,n_p\rrbracket\}
    \big\},
\end{align*}
where $n_{\mathrm{ref}}$ denotes the total number of the reference unit-cell basis functions and $n_p$ is the number of basis functions for each patch. 
%
By pulling back the macro-mapping and the deformation mapping, we can define the space of $\mathcal{C}^0$-continuous basis functions associated with the current configuration of the entire computational domain as
\begin{align*}
\begin{split}
    \spn\big\{  \varphi_{\corclem{a}} ~ | ~ {\corclem{a}}\in \llbracket1,n\rrbracket \big\}
     &= \big\{ 
        v \in \mathcal{C}^0(\Omega)
        ~|~
        \forall{s} \in \llbracket1,N_s\rrbracket, v_{|\Omega^{s}} \in \spn\{  \phi_{\corclem{a}}^{\mathrm{ref}}\circ(\bm{\Phi}\circ \mathcal{B}^{s})^{-1} ~ | ~ {\corclem{a}}\in \llbracket 1,n_{\mathrm{ref}}\rrbracket \}
    \big\},
\end{split}
\end{align*}
where $n$ denotes the total number of basis functions.
}

\modtibo{%
To represent the vector-valued displacement field in $d$ dimensions, each scalar spline basis function is associated with a $d$-dimensional vector.
\begin{defin}[Finite-dimensional approximation space]
The approximation space associated to the IGA discretization of our problem is therefore defined by
\begin{equation*}
    \bm{V}_{h}
    :=
    \left\{ 
        \bm{v} \in \bm{H}^1_0(\Omega, \partial \Omega_D) 
        ~ | ~ 
        \bm{v} \in \spn \left\{ {\varphi}_{a}\mathbf{e}_i \ | \ ~ (a,i) \in \llbracket{}1, n\rrbracket{}\times{} \llbracket{}1, d\rrbracket{}\right\} 
    \right\}. 
\end{equation*}
\end{defin}
}

\begin{rmk}
The subscript $h$ in the approximation space $\bm{V}_{h}$ refers to the micro-model spline knot-span size (i.e., the micro-model mesh size). This convention is used throughout the paper.
\end{rmk}

\newcommand{\localPhih}{\bm{\Phi}^{h,s}}

Using the notations introduced above, the global trial function representing the displacement field can be written as:
\begin{align}
\label{eq:u_phi}
\bm{u}_h(\mathbf{x}) = \sum_{a=1}^{n} \varphi_{a}(\mathbf{x})\,\mathbf{u}_{a}, \quad \forall{}\mathbf{x}\in\Omega,
\end{align}
where $\mathbf{u}_{a} \in \mathbb{R}^d$ contains the $d$ displacement coefficients associated with the $a$-th scalar basis function $\varphi_a$, representing the contributions (DOFs) in every spatial direction. Representation~\cref{eq:u_phi} is termed global because it is defined over the entire computational domain~$\Omega$. Alternatively, the global displacement field can be expressed as the assembly of local contributions associated with each cell in the lattice structure. The local displacement field $\bm{u}_h^s$, for $s\in\llbracket1,N_s\rrbracket$, is the local trial function defined over the $s$-th unit-cell $\Omega^s$ and is thus expressed using local basis functions as:
\begin{align}
    \bm{u}_h^s(\mathbf{x})
    =  
    \sum_{a=1}^{n_\mathrm{ref}}
    \varphi_a^s(\mathbf{x})\,\mathbf{u}_{a}^s,
    \quad \forall{}\mathbf{x}\in\Omega^s,
    \quad\mathrm{with}~~ 
    \varphi_a^s=\phi_{a}^{\mathrm{ref}}\circ(\localPhih{}\circ \mathcal{B}^{s})^{-1}.
    \label{eq:discr_u_s}
\end{align}
Here, $\localPhih{}$, for $s\in\llbracket1,N_s\rrbracket$, is the discrete deformation mapping associated to the $s$-th unit-cell $\Omega^s$. The functions~$\varphi_a^s$, for $a\in\llbracket1,n_\mathrm{ref}\rrbracket$, are the micro-mapping spline functions pulled-back from the current configuration domain to the reference unit cell. 

In what follows, all discrete quantities will be defined at the level of the unit-cells $\Omega^s$, thereby setting the stage for the fast assembly strategy and the efficient solver developed in \cref{sub:fast_assembly,sub:efficient_solver}, respectively. Importantly, no operator will ever be assembled or solved at the scale of the entire structure; instead, only a limited number of complete assemblies and solutions at the cell level will be required in our strategy.

\subsubsection{Residual vector: internal and external virtual work contributions}

\modtibo{%
In the Newton–Raphson procedure, the evaluation of the residual vector is required at each iteration. This vector consists of the contributions from internal and external virtual work, and will be obtained in our context by summing the local contributions associated with each cell:
\begin{equation}
    \mathbf{r}(\mathbf{u})
    = 
    \mathbf{f}^{\mathrm{int}}(\mathbf{u}) - \mathbf{f}^{\mathrm{ext}}
    =
    \sum_{s=1}^{N_s}
    \mathbf{A}^s\,(\mathbf{f}^{\mathrm{int},s}(\mathbf{u}^s) - \mathbf{f}^{\mathrm{ext},s}),
    \label{eq:residualvector}
\end{equation}
where ${\bf A}^s \in \mathbb{R}^{d{n}\,\times\,d{n}_{\mathrm{ref}}}$, for $s\in\llbracket1,N_s\rrbracket$, are the local-to-global assembly operators, mapping local DOFs to the global system. Each $\mathbf{A}_s$ is a sparse Boolean matrix containing only zeros and ones.

The components of the local load vectors are coming from the (local) virtual work functions. Specifically, the entry of the internal load vector~$\mathbf{f}^{\mathrm{int},s}$ associated with test function~$\bm{v} = \varphi_a^s\mathbf{e}_i$ is given by:
\begin{equation*}
    \mathrm{f}_{a,i}^{\mathrm{int},s}(\mathbf{u}^s)
    =
    W_\mathrm{int}^s\big(\bm{\Phi}^{h,s},\varphi_a^s\mathbf{e}_i\big)
    =
    \int_{\Omega^s} \bm{\sigma}(\bm{\Phi}^{h,s}) : \delta \bm{d}(\varphi_a^s\mathbf{e}_i) \, d{\bf x}
    =
    \int_{\Omega^s}
    {
        \sum_{j=1}^{d}
        {\sigma}_{ij}(\bm{\Phi}^{h,s}) \frac{\partial\varphi_a^s}{\partial{}x_j}
    }
    \,d{\bf x}.
\end{equation*}
In other words, the nodal sub-vector~$\mathbf{f}_a^{\mathrm{int},s}\in\mathbb{R}^d$ associated to the $a$-th local basis function is
\begin{equation}
    \mathbf{f}_a^{\mathrm{int},s}(\mathbf{u}^s)
    =
    \int_{\Omega^s}
    {
        \bm{\sigma}(\bm{\Phi}^{h,s}) \bm{\nabla}_{\mathbf{x}} \varphi_a^s
    }
    \,d{\bf x},
    \label{eq:internalloadvec}
\end{equation}
and the local vector~$\mathbf{f}^{\mathrm{int},s}$ consists in their collection for~$a = 1,\dots,n_{\mathrm{ref}}$. For a neo-Hookean material, the evaluation of the Cauchy stress tensor needed in~\cref{eq:internalloadvec} involves the computation of the deformation gradient (see~\cref{eq:sig_eq} and~\cref{eq:def_grad}). To do so, let us also express the displacement field expressed using the previously introduced spline basis functions mapped to the initial configuration:
\begin{equation*}
    \bm{u}_h^{s}(\mathbf{X})
    =
    \sum_{a=1}^{n_\mathrm{ref}}  \tilde{\varphi}_a^s(\mathbf{X})\,\mathbf{u}_{a}^{s},
    \forall{}\mathbf{X}\in\Omega_0^s, 
    \quad\mathrm{with}~~ 
    \tilde{\varphi}_a^s
    =
    {\varphi}_a^s\circ \bm{\Phi}\corclem{^{h,s}}
    = 
    \phi_{a}^\mathrm{ref}\circ(\mathcal{B}^{s})^{-1}.
\end{equation*}
With such local basis functions in hand, the deformation gradient can be expressed in indicial form as:
\begin{equation*}
    F_{ij}^{h,s}({\bf X})
    = 
    \delta_{ij}
    + 
    \sum_{a=1}^{n_\mathrm{ref}} 
    \frac{\partial \tilde{\varphi}_a^s}{\partial X_{j}}(\mathbf{X}) 
    \,{u}_{a,i}^{s},
\end{equation*}
where~${u}_{a,i}^{s} = \mathbf{u}_{a}^{s}\cdot\mathbf{e}_i$ is the $i$-th component of DOF sub-vector~$\mathbf{u}_{a}^{s}$. Finally, let us provide the expression of sub-vectors~$\mathbf{f}_a^{\mathrm{ext},s}\in\mathbb{R}^d$ which comes from the (local) external virtual work~\cref{eq:externalvirtwork} as:
\begin{equation*}
    \mathbf{f}_a^{\mathrm{ext},s}
    =
    \int_{\Omega_0^s} \bm{f}_0 \tilde{\varphi}_a^s \, d{\bf X} 
    + 
    \int_{\partial \Omega_{N,0}^s} \bm{g}_0 \tilde{\varphi}_a^s \, d{\mathbf{S}}.
\end{equation*}
Again, the local vector~$\mathbf{f}^{\mathrm{ext},s}$ consists in their collection for~$a = 1,\dots,n_{\mathrm{ref}}$.
}%

\subsubsection{Tangent operator}

\modtibo{%
The tangent operator required to form the tangent linear system at each iteration of the Newton–Raphson method is obtained by discretizing the directional derivative of the virtual work functional (see \cref{eq:direc_deriv}) using IGA. As for the residual vector~\cref{eq:residualvector}, the global tangent stiffness matrix is viewed in our context as the sum of local, cell-wise, matrices:
\begin{equation*}
    \mathbf{K}(\mathbf{u})
    = 
    \sum_{s=1}^{N_s} \mathbf{A}^s \mathbf{K}^s(\mathbf{u}^s) \mathbf{A}^{s\top},
    \quad\mathbf{K}(\mathbf{u})\in \mathbb{R}^{dn\,\times\,dn}.
\end{equation*}
Each local tangent operator is composed of two contributions, as its total continuous counterpart given in~\cref{eq:dir_der}. Specifically, a given component of the local matrices (associated to test function~$\varphi_a^s\mathbf{e}_i$ and trial function~$\varphi_b^s\mathbf{e}_j$) reads as:
\begin{equation*}
    \left[\mathbf{K}_{ab}^s(\mathbf{u}^s)\right]_{ij}
    =
    \left[\mathbf{K}_{ab}^{\mathrm{c},s}(\mathbf{u}^s)\right]_{ij}
    +
    \left[\mathbf{K}_{ab}^{\sigma,s}(\mathbf{u}^s)\right]_{ij},
\end{equation*}
where the first term is related to the constitutive behavior of the material and leads to the so-called material stiffness operator, and the second term arises from the current stress state and leads to the so-called geometrical stiffness operator. They are expressed as:
\begin{alignat*}{6}
    &
    \left[\mathbf{K}_{ab}^{\mathrm{c},s}(\mathbf{u}^s)\right]_{ij}
    &\,=\,&
    \int_{\Omega^s} \delta\mathbf{d}(\varphi_a^s\mathbf{e}_i) : \bm{\mathcal{C}}(\bm{\Phi}^{h,s}) : \bm{\varepsilon}(\varphi_b^s\mathbf{e}_j)\,d{\bf x} 
    &\,=\,&
    \int_{\Omega^s}
    \sum_{k,l=1}^d
    \frac{\partial \varphi_a^s}{\partial x_{k}}
    \mathcal{C}_{ikjl}(\bm{\Phi}^{h,s})
    \frac{\partial \varphi_b^s}{\partial x_{l}}\,d{\bf x},
    \\
    &
    \left[\mathbf{K}_{ab}^{\sigma,s}(\mathbf{u}^s)\right]_{ij}
    &\,=\,&
    \int_{\Omega^s} \bm{\sigma}(\bm{\Phi}^{h,s}) : \big(\bm{\nabla}_{\mathbf{x}}^T (\varphi_a^s\mathbf{e}_i)  \bm{\nabla}_{\mathbf{x}} (\varphi_b^s\mathbf{e}_j) \big)\,d{\bf x}
    &\,=\,&
    \int_{\Omega^s}
    \delta_{ij}
    \sum_{k,l=1}^d
    \sigma_{kl}(\bm{\Phi}^{h,s}) 
    \frac{\partial \varphi_a^s}{\partial x_{k}}
    \frac{\partial \varphi_b^s}{\partial x_{l}}\,d{\bf x}.
\end{alignat*}
Interestingly, we can define the following four-order tensor-valued function:
\begin{equation*}
    \bm{\mathcal{D}}(\bm{\Phi}^{h,s}) = \mathcal{D}_{ijkl}(\bm{\Phi}^{h,s})\,\mathbf{e}_{i}\otimes\mathbf{e}_{j}\otimes\mathbf{e}_{k}\otimes\mathbf{e}_{l},
    \quad
    \mathcal{D}_{ijkl}(\bm{\Phi}^{h,s})
    =
    \mathcal{C}_{ikjl}(\bm{\Phi}^{h,s})
    +
    \delta_{ij}
    \sigma_{kl}(\bm{\Phi}^{h,s}),
\end{equation*}
such that we can express the sub-matrices~$\mathbf{K}_{ab}^{s}\in\mathbb{R}^{d\times{}d}$ compactly as
\begin{equation*}
    \mathbf{K}_{ab}^s(\mathbf{u}^s)
    =
    \int_{\Omega^s}
    \bm{\mathcal{D}}(\bm{\Phi}^{h,s}) : 
    \big(\bm{\nabla}_{\mathbf{x}}\varphi_a^s\otimes\bm{\nabla}_{\mathbf{x}}\varphi_b^s\big)
    \,d{\bf x}.
\end{equation*}
Finally, the local matrix~$\mathbf{K}^s$ consists in their collection for~$a = 1,\dots,n_{\mathrm{ref}}$ and~$b = 1,\dots,n_{\mathrm{ref}}$.

}

\subsubsection{Gauss-Legendre quadrature}

\modtibo{%
In this work, the evaluation of the local operators entering the tangent matrix and the internal/external work vectors relies on Gauss-Legendre quadrature. In what follows, we illustrate the procedure for the tangent matrix, noting that the computation of the internal and external work vectors is entirely analogous.

The evaluation of the operators begins by pulling back the mappings~$\bm{\Psi}^{h,s} = \bm{\Phi}^{h,s}\circ\mathcal{B}^s$ so that the integrals are expressed over the reference domain~$\Omega_{\mathrm{ref}}$. It reads as:
\begin{equation}
    \mathbf{K}_{ab}^s (\mathbf{u}^s)
    =
    \int_{\Omega^s}
    \bm{\mathcal{D}}(\bm{\Phi}^{h,s}) : 
    \big(\bm{\nabla}_{\mathbf{x}}\varphi_a^s\otimes\bm{\nabla}_{\mathbf{x}}\varphi_b^s\big)
    \,d{\bf x}
    =
    \int_{\Omega_{\mathrm{ref}}}
    \bm{\mathcal{D}}^{\mathrm{ref}}(\bm{\Psi}^{h,s}) : 
    \big(\bm{\nabla_{\xi}}\phi_{a}^{\mathrm{ref}}\otimes\bm{\nabla_{\xi}}\phi_{b}^{\mathrm{ref}}\big)
    \,d\bm{\xi},
    \label{eq:stiffmatPulledback}
\end{equation}
where the components of the four-order tensor valued 
function~$\bm{\mathcal{D}}^{\mathrm{ref}}(\bm{\Psi}^{h,s}) : \Omega_{\mathrm{ref}}\to\mathbb{R}^{d\times{}d\times{}d\times{}d}$ 
are given by:
\begin{equation}
    {\mathcal{D}}^{\mathrm{ref}}_{ijkl}(\bm{\Psi}^{h,s})
    =
    \mathcal{D}_{ij\bar{k}\bar{l}}(\bm{\Phi}^{h,s})
    \left[\mathbf{J}^{-1}({\bm{\Psi}^{h,s}})\right]_{\bar{k}k}
    \left[\mathbf{J}^{-1}({\bm{\Psi}^{h,s}})\right]_{\bar{l}l}
    |\det{\mathbf{J}({\bm{\Psi}^{h,s}}})|,
    \label{eq:tensorfunctionPulledback}
\end{equation}
with Einstein's summation convention on~$\bar{k}$ and~$\bar{l}$, and $\mathbf{J}({\bm{\Psi}^{h,s}})\in\mathbb{R}^{d\times{}d}$ the Jacobian matrix of~$\bm{\Psi}^{h,s}$.

Let $\mathcal{Q}=\{(\omega_q, \bm{\xi}_q) \in \mathbb{R}\times\mathbb{R}^d  , q=1,\dots,n_g\}$ define a quadrature rule over the reference unit-cell~$\Omega_{\mathrm{ref}}$. The local tangent stiffness matrices can be built as:
\begin{equation}%
    \mathbf{K}_{ab}^s(\mathbf{u}^s)
    \approx
    \sum_{q=1}^{n_g}
    \omega_q\,
    \bm{\mathcal{K}}_{ab}^s(\bm{\xi}_q),
    \quad\mathrm{with}\quad
    \bm{\mathcal{K}}_{ab}^s(\bm{\xi}_q)
    =
    \bm{\mathcal{D}}^{\mathrm{ref}}(\bm{\Psi}^{h,s})(\bm{\xi}_q) : 
    \big(\bm{\nabla_{\xi}}\phi_{a}^{\mathrm{ref}}(\bm{\xi}_q)\otimes\bm{\nabla_{\xi}}\phi_{b}^{\mathrm{ref}}(\bm{\xi}_q)\big),
    \label{eq:fullquadraturerule}
\end{equation}%
and $a = 1,\dots,n_{\mathrm{ref}}$ and~$b = 1,\dots,n_{\mathrm{ref}}$. The generation of the quadrature rule is done classically in an element-by-element fashion, using Gauss-Legendre quadrature rule with~$(p+1)^d$ points per element.

}

\subsection{Fast-assembly with a reduced basis strategy}
\label{sub:fast_assembly}

\begin{figure}
    \centering
    \includegraphics[width=0.8\linewidth]{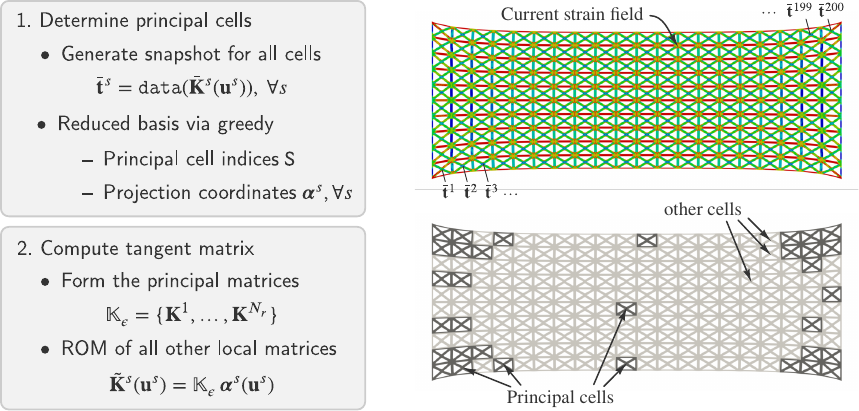}%
	\caption{RB strategy to compute local tangent operators.}
	\label{fig:RBstrategy}
\end{figure}

Due to the inherent use of high-degree basis functions in IGA, the computational cost, both in terms of time and memory, of assembling and solving the global tangent stiffness matrix at each Newton–Raphson iteration can be substantial. The goal in the following is therefore to propose an efficient assembly and solver strategy that exploits the repetitive nature of lattice geometries, in the same spirit as previous contributions~\cite{Hirschler:2023aa,GUILLET2025114136}, but extended here to nonlinear analysis.

The core idea of the proposed fast assembly strategy is to exploit the inherent repetitiveness of the lattice cells in order to reduce the number of local (cell-level) computations. In other words, since lattice structures are generated by repeatedly applying a macroscopic mapping to a common reference unit cell, strong structural similarities naturally arise in the local discrete operators. These similarities can be leveraged to reduce computational cost and memory requirements.
The proposed strategy, detailed hereafter, is illustrated in \cref{fig:RBstrategy}. First, a reduced basis (RB) is constructed from a limited set of so-called principal cells. Then, each local operator can then be accurately approximated within this RB space, allowing in particular the tangent matrix to be efficiently formed and stored at reduced cost. 

\subsubsection{Reduced basis approximation}
\label{sec:rb}

Let us introduce the following matrix manifold, defined as the collection of all local tangent matrices associated with the cells of the lattice structure:
\begin{equation*}
    \mathbb{K} = \left\{\mathbf{K}^s(\mathbf{u}^s),  s = 1,\dots, N_s\right\}.
\end{equation*}
Thanks to the repetitive nature of the lattice cells, strong structural similarities are expected among the local tangent matrices. Consequently, the effective dimension of this manifold is assumed to be significantly lower than the total number of cells~$N_s$. To exploit this property, we aim to approximate $\mathbb{K}$ by a lower-dimensional manifold, associated with a subset of cells which can be defined as the principal cells. Specifically, to achieve this we denote the list of indices of the considered principal cells by~$\mathsf{S} = \{ s_1, \dots, s_{N_r} \}$ and we introduce a reduced basis as:
\begin{equation}
    \mathbb{K}_{\epsilon}
    = 
    \{\mathbf{K}^1, \dots, \mathbf{K}^{N_r}\},
    \quad\mathrm{with}\quad
    \mathbf{K}^r
    =
    \mathbf{K}^{s_r}(\mathbf{u}^{s_r}).
    \label{eq:reducedbasis}
\end{equation}
The associated reduced space is then defined by
\begin{equation*}
    \mathscr{K}_{\epsilon} := \spn(\mathbb{K}_{\epsilon}),
\end{equation*}
such that:
\begin{align}
    \forall{}\mathbf{K}^s(\mathbf{u}^s) \in \mathbb{K}, 
    \quad \exists \tilde{\mathbf{K}}^s(\mathbf{u}^s) \in \mathscr{K}_{\epsilon}, 
    ~~\mathrm{satisfying}~~ 
    \| \mathbf{K}^s(\mathbf{u}^s) - \tilde{\mathbf{K}}^s(\mathbf{u}^s) \|_{\star} < \epsilon,
    \label{eq:tol}
\end{align}
where $\|\cdot\|_{\star}$ denotes an appropriate matrix norm (e.g., Frobenius or operator norm). Finally, an approximated local tangent matrix~$\tilde{\mathbf{K}}^s\in\mathscr{K}_{\epsilon}$ is expressed as a linear combination of the principal tangent matrices:
\begin{equation}
    \tilde{\mathbf{K}}^s(\mathbf{u}^s)
    = 
    \sum_{r=1}^{N_r} 
    \alpha_r^s(\mathbf{u}^s)\,\mathbf{K}^{r}
    =
    \mathbb{K}_{\epsilon}\,
    \bm{\alpha}^s(\mathbf{u}^s)
    ,
    \label{eq:romstiffnessmatrix}
\end{equation}
where $\bm{\alpha}^s(\mathbf{u}^s) \in \mathbb{R}^{N_r}$ are the RB coefficients.

The key advantage of this approach is that only the $N_r$ principal local tangent matrices need to be explicitly assembled, and then the entire set of $N_s$ local tangent matrices can be represented using the $N_s \times N_r$ scalar coefficients~$\alpha_r^s$. Moreover, as discussed in \cref{sub:efficient_solver}, only these $N_r$ principal tangent matrices are actually \corclem{used} within the developed inexact domain decomposition solver. The later use of such a dedicated solver also explains why the reduced basis is built based on the selection of principal cells (remember~\cref{eq:reducedbasis}), and not through modes as one would get from the proper orthogonal decomposition for instance.

\begin{rmk}
It is also worth emphasizing that the proposed RB-based operator approximation becomes particularly appealing in the context of a Newton–Raphson algorithm. Indeed, the operators only need to be approximated with sufficient accuracy to ensure the overall convergence of the nonlinear algorithm, rather than to precisely solve each individual tangent linear system. In other words, since these approximations are embedded within a Newton scheme, the required level of accuracy can be relaxed compared to what would be needed for solving a single linear problem.
\end{rmk}

\subsubsection{Construction of the reduced basis}
\label{sec:constr_rb}

In practice, the reduced basis~$\mathbb{K}_{\epsilon}$ introduced in~\cref{eq:reducedbasis} is constructed using a greedy algorithm. The latter can be interpreted as a truncated Gram-Schmidt orthogonalization. It operates on a collection of so-called snapshots as standard algorithms for RB construction. Several options of snapshots seem to be possible in our context, inspired by the Empirical Interpolation Method and its discrete variant~\cite{maday2013generalized,quarteroni15,hesthaven16}. Here, the snapshots need to be computationally cheap as it is done online, at each Newton iteration. 

To do so, we apply an almost non-intrusive approach, in a M-DEIM-fashion~\cite{negri2015efficient}. A snapshot is made of the non-zero components of the local tangent stiffness matrix built using a Gauss-Legendre quadrature rule with less points per direction (i.e., a naive reduced quadrature rule) denoted here by~$\bar{\mathcal{Q}} = \{(\bar{\omega}_q, \bar{\bm{\xi}}_q) \in \mathbb{R}\times\mathbb{R}^d, q=1,\dots,\bar{n}_g\}$. Thus, a snapshot~$\bar{\mathbf{t}}^s$ reads as
\begin{equation}
    \bar{\mathbf{t}}^s = \mathtt{data}(\bar{\mathbf{K}}^{s}(\mathbf{u}^s))
    \quad\mathrm{with}\quad
    \bar{\mathbf{K}}^{s}(\mathbf{u}^s)
    = 
    \sum_{q=1}^{\bar{n}_g}
    \bar{\omega}_q\,
    \bm{\mathcal{K}}^s(\bar{\bm{\xi}}_q).
    \label{eq:snapshots}
\end{equation}
As will be shown in the numerical results \cref{sec:res_num}, the number of reduced Gauss points $\bar{n}_g$ per cell is chosen to be significantly smaller than the number of integration points $n_g$ required for an accurate evaluation of the local tangent operators (see~\cref{eq:fullquadraturerule}). In practice, using such a reduced quadrature rule at this stage proves sufficient, since the goal is only to identify the principal cells, rather than to compute a good approximation of the operators themselves. Also, it is weakly intrusive: one simply need to change the order of the quadrature rule during the matrix formation which is a regular input parameter of FE/IGA codes. 

\begin{rmk}
A more intrusive yet interesting choice of snapshots would consist in the evaluation of the fourth-order tensor-valued function~$\bm{\mathcal{D}}^{\mathrm{ref}}(\bm{\Psi}^{h,s})$ as defined in~\cref{eq:tensorfunctionPulledback} at given points~$\bar{\bm{\xi}}_q\in\Omega^{\mathrm{ref}}, q=1,\dots,\bar{n}_g$. Such a choice would read as
\begin{equation*}
    \bar{\mathbf{t}}^s 
    = 
    \mathtt{vec}(\{ 
        \bm{\mathcal{D}}^{\mathrm{ref}}(\bm{\Psi}^{h,s})(\bar{\bm{\xi}}_q),~q = 1,\dots,\bar{n}_g
    \}).
\end{equation*}
Indeed, based on~\cref{eq:stiffmatPulledback}, one can see that all the cell-to-cell differences are concentrated in the functions~$\bm{\mathcal{D}}^{\mathrm{ref}}(\bm{\Psi}^{h,s})$. The rest of the integrand only involves quantities associated with the reference unit cell (i.e., being identically shared by all the cells) which do not play a role in the mechanical differences between the cells. Finally, one could also think of generating the snapshots by computing a limited number of entries of the ``true'' local tangent stiffness matrix
\begin{equation*}
    \bar{\mathbf{t}}^s 
    = \mathtt{some\_data}(\mathbf{K}^{s}(\mathbf{u}^s)),
\end{equation*}
in the M-DEIM spirit.
\end{rmk}

Once all the snapshots are generated via~\cref{eq:snapshots}, the snapshots matrix~$\mathbf{T} = \begin{pmatrix} \bar{\mathbf{t}}^1 & \dots & \bar{\mathbf{t}}^{N_s} \end{pmatrix}, \mathbf{T}\in\mathbb{R}^{\bar{n}_{nz}\times{}N_s}$ is fed into a greedy algorithm. The algorithmic details are described in~\cref{alg:1}. In a nutshell, based on a given tolerance~$\epsilon$, it returns (among others) the principal cell indices~$\mathsf{S}$, the associated reduced basis~$\bar{\mathbb{t}}_{\epsilon} = (~\bar{\mathbf{t}}^{s_1}~\cdots~\bar{\mathbf{t}}^{s_{N_r}}~)\in\mathbb{R}^{\bar{n}_{nz}\times{}N_r}$, and the projection coordinates~$\bm{\alpha}^s$, i.e.,
\begin{equation*}
    \left(
        \mathsf{S}, 
        \bar{\mathbb{t}}_{\epsilon},
        \{\bm{\alpha}^s(\mathbf{u}^s),~s=1,\dots,N_s \}
    \right)
    =
    \mathtt{greedy}(\mathbf{T}, \epsilon).
\end{equation*}
Based on~$\mathsf{S}$, one can generate the final reduced basis~$\mathbb{K}_{\epsilon}$ via~\eqref{eq:reducedbasis}, and based on the~$\bm{\alpha}^s(\mathbf{u}^s)$, one can get the (well) approximated local tangent stiffness matrices~$\tilde{\mathbf{K}}^s(\mathbf{u}^s)$ via~\cref{eq:romstiffnessmatrix}. To summarize, it generates a reduced order modelling of the local operators as
\begin{equation*}
    \mathbf{K}^s(\mathbf{u}^s)
    ~\approx~
    \tilde{\mathbf{K}}^s(\mathbf{u}^s)
    =
    \mathbb{K}_{\epsilon} 
    \,\bm{\alpha}^s(\mathbf{u}^s),
    \quad\mathrm{where}\quad
    \bm{\alpha}^s(\mathbf{u}^s)
    =
    \big[\bar{\mathbb{t}}_{\epsilon}^\top\bar{\mathbb{t}}_{\epsilon}^{}\big]^{-1}\bar{\mathbb{t}}_{\epsilon}^\top\bar{\mathbf{t}}^s(\mathbf{u}^s)
    .
\end{equation*}
Unlike the sub-integrated matrices~$\bar{\mathbf{K}}^s(\mathbf{u}^s)$ involved during the snapshots generation, the ROM matrices~$\tilde{\mathbf{K}}^s(\mathbf{u}^s)$ are suitable for analysis.

\algrenewcommand\algorithmicrequire{\textbf{Input:}}
\algrenewcommand\algorithmicensure{\textbf{Output:}}
\begin{algorithm}[t]
\caption{Greedy algorithm for the identification of the principal cells.}
\label{alg:1}
\begin{algorithmic}[1]
\Require Snapshots for all cells: $\bar{\mathbf{t}}^s$
         \Statex Tolerance of the reduced basis: $\epsilon>0$

\For{$s = 1$ \textbf{to} $N_s$}
    \State $\bm{\Delta}_0^{s} \leftarrow {\bar{\mathbf{t}}^s}/{\|\bar{\mathbf{t}}^s\|_2}$
    \Comment{Normalization}
\EndFor

\State $i \leftarrow 0$
\While{$\max\limits_{s=1,\ldots,N_s} \|\bm{\Delta}_i^{s}\|_\infty > \epsilon$}
    \State $i \leftarrow i + 1$
    \State $s_i \leftarrow \arg\max\limits_{s=1,\ldots,N_s} \|\bm{\Delta}_{i-1}^{s}\|_\infty$
    \Comment{Find current principal unit cell}
    \State $\bm{\zeta}_i \leftarrow \bm{\Delta}_{i-1}^{s_i} / \|\bm{\Delta}_{i-1}^{s_i}\|_2$
    \Comment{Add new vector basis}
    \For{$s = 1$ \textbf{to} $N_s$}
        \State ${\beta}_i^{s} \leftarrow \bm{\Delta}_{i-1}^{s} \cdot \bm{\zeta}_i$
        \Comment{Compute additional coordinates}
        \State $\bm{\Delta}_i^{s} \leftarrow \bm{\Delta}_{i-1}^{s} - {\beta}_i^{s} \bm{\zeta}_i$
        \Comment{Update residuals}
    \EndFor
\EndWhile

\State $N_r \leftarrow i$

\Ensure
    Orthonormal basis: $\mathcal{Z} = \{\bm{\zeta}_i \mid i\in \llbracket1,N_r\rrbracket \}$
    \Statex Principal unit-cell indices: $\mathsf{S} = \{s_i \mid i\in \llbracket1,N_r\rrbracket\}$
    \Statex Certified affine decompositions: $\tilde{\mathbf{t}}^{s} = \|\bar{\mathbf{t}}^s\|_2 \mathcal{Z} \bm{\beta}^{s} \quad \forall s \in \llbracket1,N_s \rrbracket$
\end{algorithmic}
\end{algorithm}

\begin{figure}
    \centering
    \begin{minipage}[t]{0.26\textwidth}\centering
        \includegraphics[width=\textwidth]{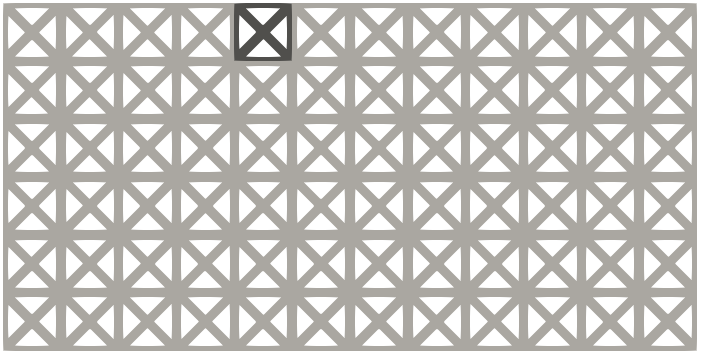}
        \centering {Linear: 1/64}
        \end{minipage}
    \begin{minipage}[t]{0.26\textwidth}\centering
        \includegraphics[width=\textwidth]{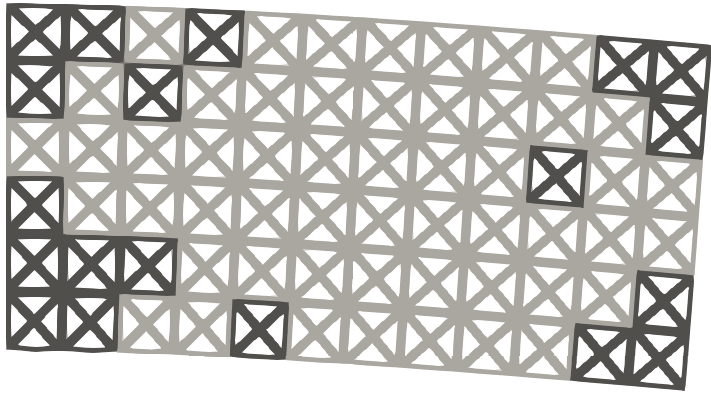}
        \centering {Nonlinear: 19/64}
        \end{minipage}
         \begin{minipage}[t]{0.26\textwidth}\centering
        \includegraphics[width=\textwidth]{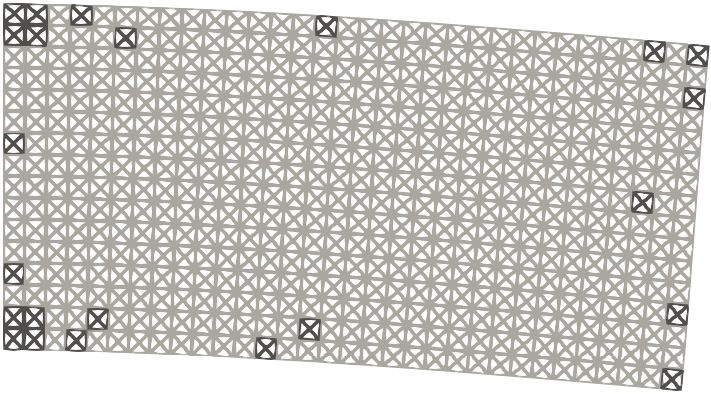}
        \centering {Nonlinear: 23/512}
        \end{minipage}
  \begin{minipage}[t]{0.34\textwidth}\centering
        \includegraphics[width=\textwidth]{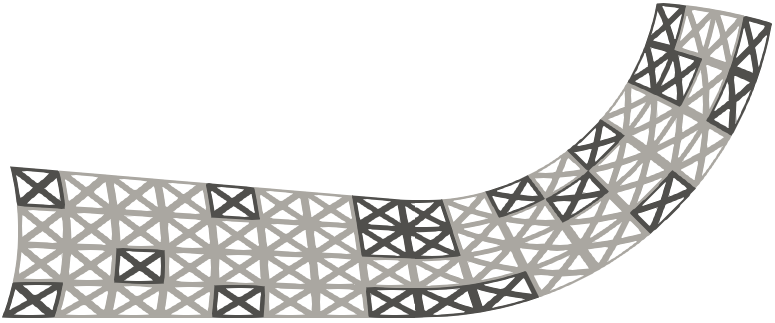}
        \centering {Linear: 21/64}
        \end{minipage}
    \begin{minipage}[t]{0.34\textwidth}\centering
        \includegraphics[width=\textwidth]{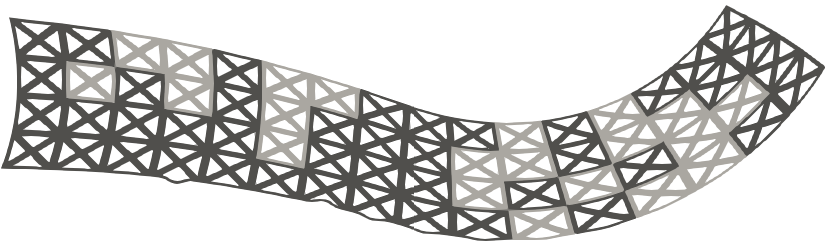}
        \centering {Nonlinear: 43/64}
        \end{minipage}
	\caption{Examples of identified principal cells for (top) a rectangular lattice and (bottom) a brake-pedal lattice, both based on a cross-truss unit cell.
The reduced-basis tolerance is set to $\epsilon = 3 \times 10^{-4}$. Results are shown for both the linear (small-deformation) and nonlinear (hyperelastic large-deformation) cases. As expected, increasing the complexity of the macroscopic mapping or the deformation level leads to a larger number of principal cells. Conversely, as the number of cells increases, the likelihood of encountering mechanically similar ones also increases, thereby enhancing the computational savings achieved by the proposed strategy.}
	\label{fig:3}
\end{figure}

Ultimately, this RB-based approach is applied at each Newton iteration. It is important to note that the “mechanical” proximity between cells, which the method exploits, also depends on the deformation experienced by the cells throughout the nonlinear simulation. The proposed strategy is all the more efficient as the number of principal cells $N_r$ required at each Newton iteration remains small. Examples of identified principal cells, obtained with a tolerance of $\epsilon = 10^{-3}$ in \cref{eq:tol}, are shown in \cref{fig:3} for two different macroscopic mappings, various numbers of cells, and distinct deformation states. Naturally, the more complex the macroscopic mapping and the induced deformations, the larger the number of principal cells needed. Conversely, as the number of cells increases, so does the likelihood of finding mechanically similar ones, which further enhances the computational savings achieved by the proposed strategy. More detailed studies are presented in \cref{sec:res_num}.

\begin{rmk}
    One could consider other hyperreduction techniques for a fast formation of the tangent operators, as for instance \cite{farhatDimensionalReductionNonlinear2014, hernandez2017dimensional, bhattacharyyaHyperReductionTechniquesEfficient2025}. However, their direct compatibility with the inexact \mbox{FETI-DP} preconditioner described hereafter does not appear straightforward to us.
\end{rmk}

\subsection{Inexact FETI-DP solver}
\label{sub:efficient_solver}

\subsubsection{Principle}

In this section we introduce a method to efficiently solve the linear systems in \cref{eq:linsys}, resulting from the discretization of the linearized equilibrium equations at each Newton iteration. Because of memory limitations, direct solvers are often impractical for large size linear systems so that we opt for an iterative solver. The idea is again to benefit from the repetitiveness of the unit cells in the structure by re-using the RB strategy introduced in the previous section at the solver stage. 
We consider a non-overlapping substructuring method based on the natural domain decomposition of the lattice structure into unit cells:  the ROM-based inexact FETI-DP method introduced in~\cite{Hirschler:2023aa}. This method involves three main components:
\begin{enumerate} [label=(\roman*)]
\item In FETI-DP methods, the continuity constraints on the displacement at the unit-cell corners are maintained throughout the iterative process, naturally leading to a coarse problem at each iteration, while other constraints are enforced using Lagrange multipliers.
\item By following the framework of inexact FETI-DP algorithms \cite{klawonn_2006}, we avoid the resolution of numerous local systems. This is achieved by iterating on the initial complete saddle-point problem and designing a dedicated block preconditioner consistent with the problem.
\item The principal unit cells identified within the RB strategy (see \cref{sec:rb}) are used in this block preconditioner to build reduced bases for efficiently approximating the numerous local systems.
\end{enumerate}

Let us briefly introduce such an inexact FETI-DP method and how the RB strategy is used within it. More details on the method can be found in \cite{Hirschler:2023aa,GUILLET2025114136}. The starting point is the following saddle-point problem:
\modtibo{%
\begin{equation*}
    \begin{pmatrix}
        \mathbf{K}(\mathbf{u}) & \mathbf{B}^\top \\
        \mathbf{B} & \mathbf{0}
    \end{pmatrix}
    \begin{pmatrix}
        \Delta\mathbf{u}  \\
        \bm{\lambda} 
    \end{pmatrix}
    =
    \begin{pmatrix}
        -\mathbf{r}(\mathbf{u}) \\
        \mathbf{d} 
    \end{pmatrix},
    \quad\mathrm{with}\quad  
    \mathbf{K}
    =
    \begin{pmatrix}
        \mathbf{K}_{RR} & \mathbf{K}_{RP} \\
        \mathbf{K}_{RP}^\top & \mathbf{K}_{PP}  
    \end{pmatrix}, 
    \quad
    \Delta\mathbf{u} =
    \begin{pmatrix}
        \Delta\mathbf{u}_{R}  \\
        \Delta\mathbf{u}_{P} 
    \end{pmatrix} 
    \quad\mathrm{and}\quad
    \mathbf{r} 
    =
    \begin{pmatrix}
        \mathbf{r}_{R}  \\
        \mathbf{r}_{P} 
    \end{pmatrix},
\end{equation*}
}%
where $\mathbf{K}(\mathbf{u})$ is the global tangent matrix and $\mathbf{r}(\mathbf{u})$ is the discrete residual vector. 
The displacement DOFs have been partitioned into primal DOFs (at cell corners, denoted by subscript $P$), which are assembled globally, and remaining DOFs (denoted by subscript $R$, which include interior DOFs, and dual DOFs at cell interfaces, the latter being enforced through Lagrange multipliers). $\bm{\lambda}\in\mathbb{R}^L$ are these Lagrange multipliers, where $L$ is the total number of discrete continuity equations imposed weakly (encapsulated in $\mathbf{d}$). $\mathbf{B}$ is the coupling matrix, associated to the cell-to-cell interface continuity conditions ($\mathbf{B}$ vanishes for the primal and interior DOFs). The idea behind inexact FETI-DP methods is to iterate on this saddle-point problem using a suitable block preconditioner. 
\modtibo{%
Following~\cite{Hirschler:2023aa}, our block-preconditioner reads as:
\begin{equation}
    \tilde{\mathbf{A}}_{\mathtt{iFETIdp}}^{-1}
    =
    \begin{bmatrix}
        \mathbf{I} & -\tilde{\mathbf{U}}^\top \\
        & \mathbf{I}
    \end{bmatrix}
    \begin{bmatrix}
        \tilde{\mathbf{K}}^{-\mathtt{rom}} \\
        & \tilde{\mathbf{F}}^{-\mathtt{pcg}}
    \end{bmatrix}
    \begin{bmatrix}
        \mathbf{I} & \\
        -\tilde{\mathbf{U}} & \mathbf{I}
    \end{bmatrix},
    \label{eq:ifetidp_precond}
\end{equation}
where local reduced order models are involved in the sub-blocks~$\tilde{\mathbf{U}}, \tilde{\mathbf{K}}, \tilde{\mathbf{F}}$. These ROMs approximate the DD local operators and are generated based on the already identified principal cells (during the fast-assembly strategy, see~\cref{sub:fast_assembly}). For instance, in the case of the local dual Schur complement (denoted here as~$\mathbb{F}_{dd}$), we generate a RB as:
\begin{equation*}
    \mathbb{F}_{dd} = \{ \mathbf{F}_{dd}^1, \dots, \mathbf{F}_{dd}^{N_r}\},
\end{equation*}
where~$\mathbf{F}_{dd}^r$ is the dual Schur complement associated to principal cell~$s_r\in\mathsf{S}$. We then build projection-based reduced-order models of the following form
\begin{equation*}
    \mathbf{F}_{dd}^s(\mathbf{u}^s)
    ~\approx~
    \tilde{\mathbf{F}}_ {dd}^s(\mathbf{u}^s)
    =
    \mathbb{F}_{dd}
    \,\bm{\delta}^s(\mathbf{u}^s),
\end{equation*}
where the projection coefficients~$\bm{\delta}^s(\mathbf{u}^s)\in\mathbb{R}^{N_r}$ are obtained classically by solving small $N_r$-by-$N_r$ linear systems~\cite[see][Section 3.1.2]{Hirschler:2023aa}.
}%

\begin{figure}
    \centering 
    \includegraphics[width=0.86\textwidth]{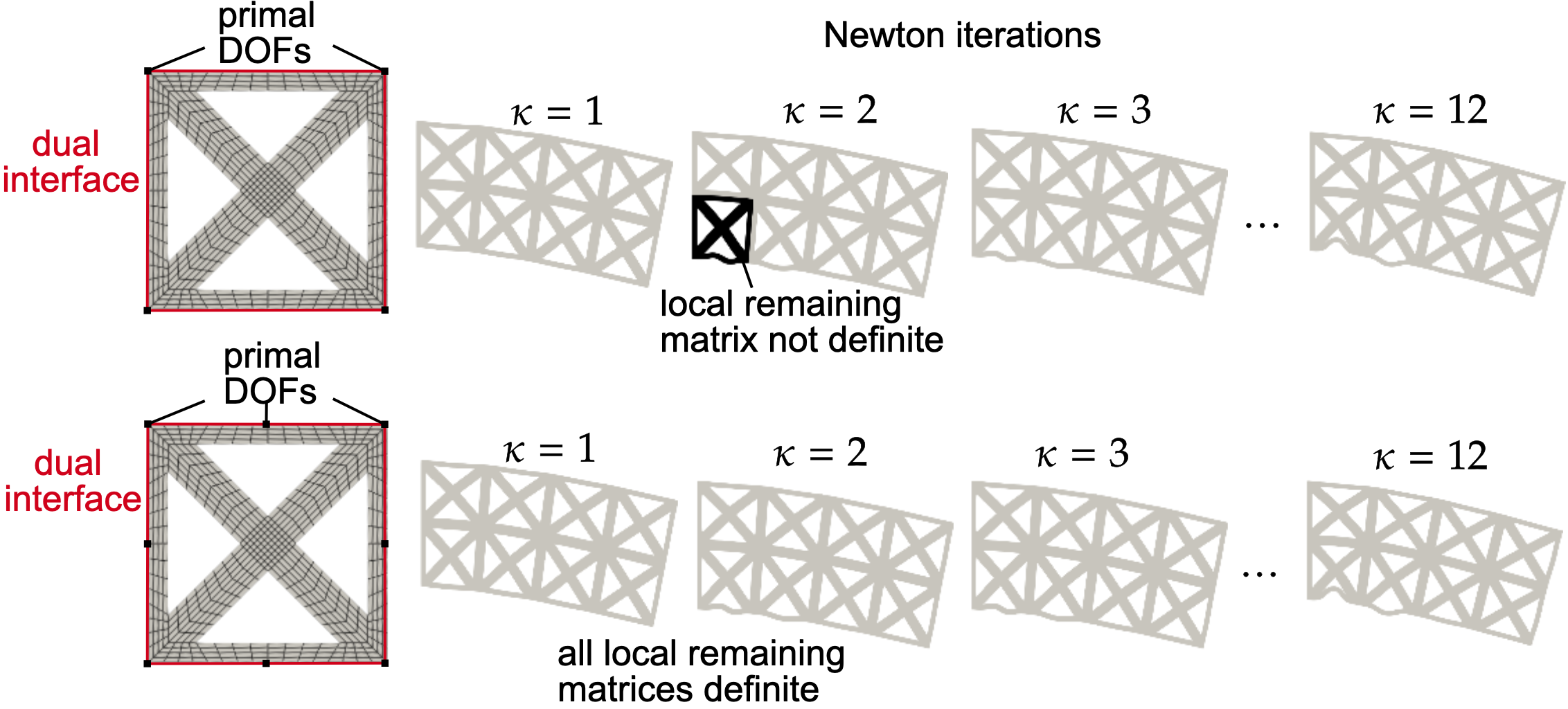}
    \caption{Non definiteness of a local remaining matrix for a deformation with buckling of bottom-hand unit-cell edge under the action of a vertical downward force applied to the right face of the structure. All matrices become definite when adding one primal DOF on each edge of the unit cell.} 
    \label{fig:4}
\end{figure}

\subsubsection{Enhanced coarse problem}
For three-dimensional problems or problems with large deformations, constructing the coarse problem by selecting only DOFs at unit-cell corners may not be sufficient to ensure invertibility of local tangent matrices or to ensure a small condition number and achieve convergence within a reasonable number of iterations for the interface problem~\modtibo{(i.e., the application of the sub-preconditioner~$\tilde{\mathbf{F}}^{-\mathtt{pcg}}$ in~\cref{eq:ifetidp_precond})}. 

First, in order to ensure a small condition number of the interface problem, we adopt the method introduced in~\cite{Farhat_2000}, which augments the coarse problem by adding edge averages. Note at this stage that we nevertheless do not add first-order moments, as we are focusing on homogeneous problems without coefficient jumps across cell boundaries. The enforcement of edge averaging is done by introducing a set of additional Lagrange multipliers for all the edges at the interfaces between cells, denoted as $\bm{\mu} \in \mathbb{R}^{E}$, where $E$ is the number of such edges. The primal variables are unchanged by this modification and remain associated with the cell corners.
We also introduce a matrix $\mathbf{Q} \in \mathbb{R}^{L \times E}$ corresponding to the translation rigid body modes of each edge in the $x$-, $y$-, and $z$-directions, see~\cite{Farhat_2000} for details. The resulting saddle-point problem for the augmented coarse problem is then formulated as follows:
\begin{align*}
    \begin{pmatrix}
        \mathbf{K}(\mathbf{u}) & \mathbf{B}^\top &\mathbf{B}^\top\mathbf{Q}  \\
        \mathbf{B} & \mathbf{0} & \mathbf{0}  \\
        \mathbf{Q}^\top\mathbf{B} & \mathbf{0}& \mathbf{0} \\
    \end{pmatrix}
    \begin{pmatrix}
        \Delta\mathbf{u}  \\
        \bm{\lambda} \\
        \bm{\mu}
    \end{pmatrix}
    =
    \begin{pmatrix}
        -\mathbf{r}(\mathbf{u}) \\
        \mathbf{d} \\
        \mathbf{Q}^\top\mathbf{d}
    \end{pmatrix}.
\end{align*} 
By eliminating the remaining primal DOFs, $\mathbf{u}_R$ and $\mathbf{u}_P$, along with the additional Lagrange multipliers $\bm{\mu}$, an interface problem including the augmented coarse problem is constructed. 

\modtibo{Second, we increase the size of the coarse problem by selecting new primal variables on the edges (or face in three dimensions) of the unit cells, during the Newton iteration.} This modification is based on the observation that local remaining (i.e. non-primal) tangent matrices associated to unit cells that undergo large deformation may not be definite and thus invertible when selecting only the corner DOFs as primal variables. A typical scenario when this happens is for example, when the struts of a crossed-shaped unit cell buckle under the action of the load applied. This situation is illustrated in \cref{fig:4}, where local buckling occurs at the bottom-hand edge of one unit cell. This phenomenon leads to a non-definite local remaining matrix of cell~$s=0$, i.e. $\mathbf{K}_{rr}^0(\mathbf{u}^0)$, at iteration $\currentiter=2$. To restore definiteness and ensure invertibility, we enrich the set of primal variables by adding one DOF on each edge of the unit cell. Importantly, the loss of definiteness arises only at iteration $\currentiter=2$ and does not persist in subsequent iterations.
In this work, we initialize the simulations with four primal variables located at the corners of the unit cell. Whenever a local matrix becomes non-definite, we systematically add one primal variable on each edge of all unit cells, and keep this configuration for the remainder of the simulation. If non-definiteness occurs again, additional edge variables are introduced in the same way until all local matrices are definite. This global enrichment strategy is motivated by the fact that, for the RB approach in FETI-DP to remain optimal, all unit cells must share an identical partitioning of primal, dual, and remaining DOFs. Consequently, even in cases such as \cref{fig:4}, where only a few unit cells are non-definite, primal variables are consistently added to every unit cell.

\begin{figure}[pos=h]
    \centering 
    \includegraphics[width=\linewidth]{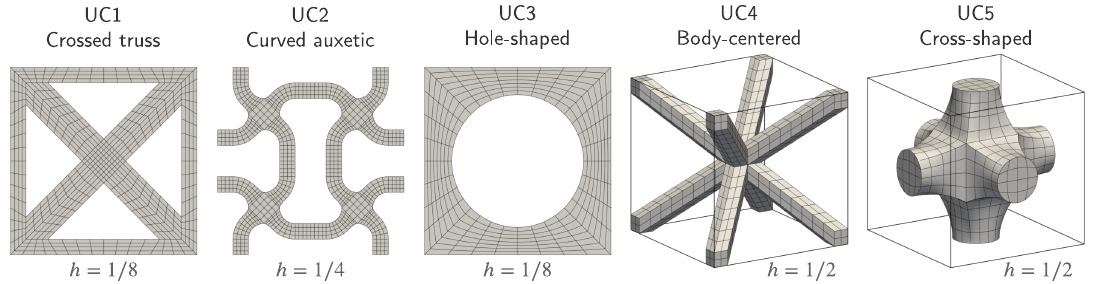}%
    \caption{Different analysis-suitable models (obtained by applying $k$-refinement) of the unit cells (UC) considered. Characteristic mesh size is defined as~$h = 1/n_{e}$, where $n_e$ is approximately the number of element per direction per patch of the underlying model.}
    \label{figparam:1}
\end{figure}

\begin{figure}[pos=t]
    \centering
    \includegraphics[width=\linewidth]{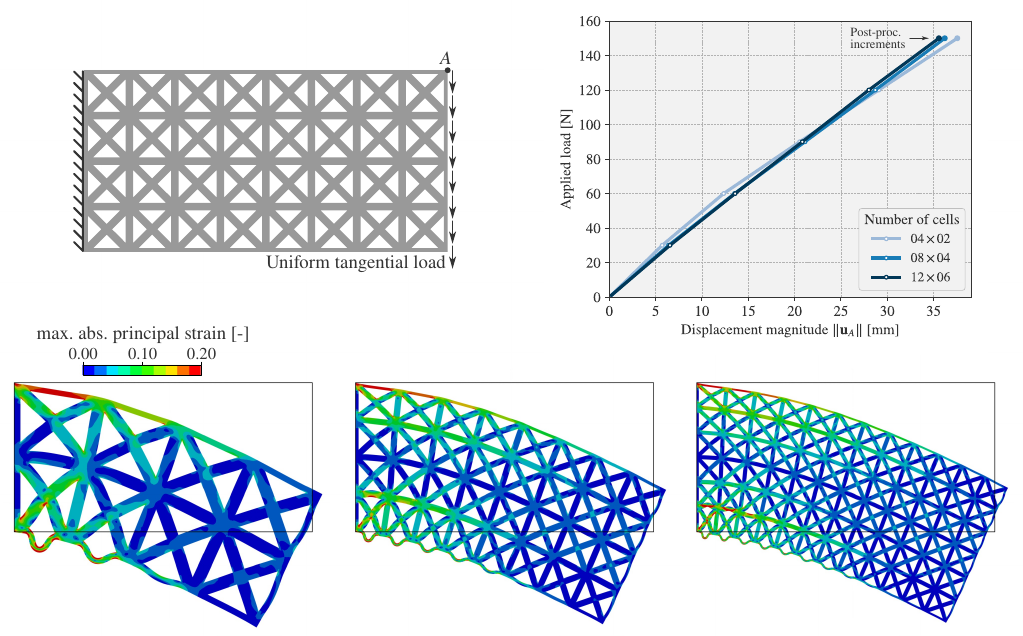}%
    \caption{2D planar map with crossed truss UC1 unit cells: a vertical load of $g_y$ is applied to the right edge of the structure. The lattice structure consists of $N_s=32$ unit cells with $p=3$, $h=1/8$.}
    \label{fignum:1}
\end{figure}

\section{Numerical results}
\label{sec:res_num}
In this section, we investigate several applications. For each case, the computational domain, i.e., the lattice structure, is 
constructed from a reference (microscopic) unit cell (UC) combined with a macroscopic geometry, as detailed in \cref{sec:1.2}. Five distinct reference UCs are employed (see \cref{figparam:1}) each exhibiting specific geometric and topological characteristics.

\begin{itemize} 
\item[(UC1)] Crossed truss unit cell: The first micro-geometry consists of a two-dimensional hollow square linked by a cross and is composed of straight struts arranged in a diagonally intersecting truss pattern. The UC is initially constructed from 24 linear Bézier patches. 
\item[(UC2)]  Curved auxetic unit cell: The micro-geometry is a two-dimensional auxetic lattice pattern, characterized by smoothly curved struts oriented to form a re-entrant geometry. The UC is initially modeled with a total of 20 matching patches, consisting of 16 quadratic NURBS patches and 4 linear patches. The curvature of the elements induces auxetic behavior, enabling, in particular, lateral expansion under tensile axial loading. 
\item[(UC3)]  Hole-shaped unit cell: The micro-geometry is a two-dimensional hollow square with a circular hole, initially constructed with 4 quadratic patches.

\item[(UC4)] Body-centered unit cell: The micro-geometry is modeled as a three-dimensional body-centered cubic (BCC) cell using, initially, 32 linear patches. The UC features struts connecting the center of the cubic domain to each corner node.

\item[(UC5)] Cross-shaped unit cell: The micro-geometry is a three-dimensional cross-shaped pattern, initially constructed using 7 cubic patches. The UC comprises orthogonal struts intersecting at the center of the cube, aligned with the principal Cartesian axes.
\end{itemize}
All of these reference unit cells are then further refined using $k$-refinement, that is, by first elevating the degree and then applying knot insertion. The resulting models are analysis-suitable, ensuring matching interfaces between patches and avoiding trimmed geometries. Examples of such analysis models are given on \cref{figparam:1}. We consider a hyperelastic neo-Hookean material as the constitutive model (see \cref{sec:const_eq}), with Young's modulus~$E=\qty{500}{\mega\Pa}$ and Poisson's ratio~$\nu=0.40$

\begin{rmk}
    Unit-cells UC2 and UC5 do not have DOF at corners of the underlying unit square/cube. Primal DOF are then selected based on closest control points to the corners, as done in~\citet{Hirschler:2023aa}.
\end{rmk}

In \cref{sec:num:1}, we analyze the mechanical response of the lattice structures, using our efficient, RB-based approach. In contrast, in \cref{sec:num:2} we focus on solver performance. For this purpose, we compare our RB method---based on an RB strategy for the efficient computation of local operators and the solution of the linear system using an inexact FETI-DP solver--- to a method using standard computation of the tangent operator and LU factorization. 
We restrict our investigation to a sequential implementation, with simulations performed on a single CPU core. All the numerical experiments are computed on a  3,5 GHz M2 (Apple) processor with 24 GB RAM. The code is implemented in Python using petsc4py~\cite{petsc-user-ref} as the linear algebra backend, together with the libraries YETI~\cite{Duval_2023} and IGAlattice.

\subsection{Structural analysis}
\label{sec:num:1}

In this section, we investigate the structural behavior of several lattice structures by solving the equilibrium equations under prescribed external forces or imposed displacements. In the following figures, we post-process the logarithmic strains (or Hencky strains).


\begin{figure}[pos=t]
    \centering
    \includegraphics[width=0.9\linewidth]{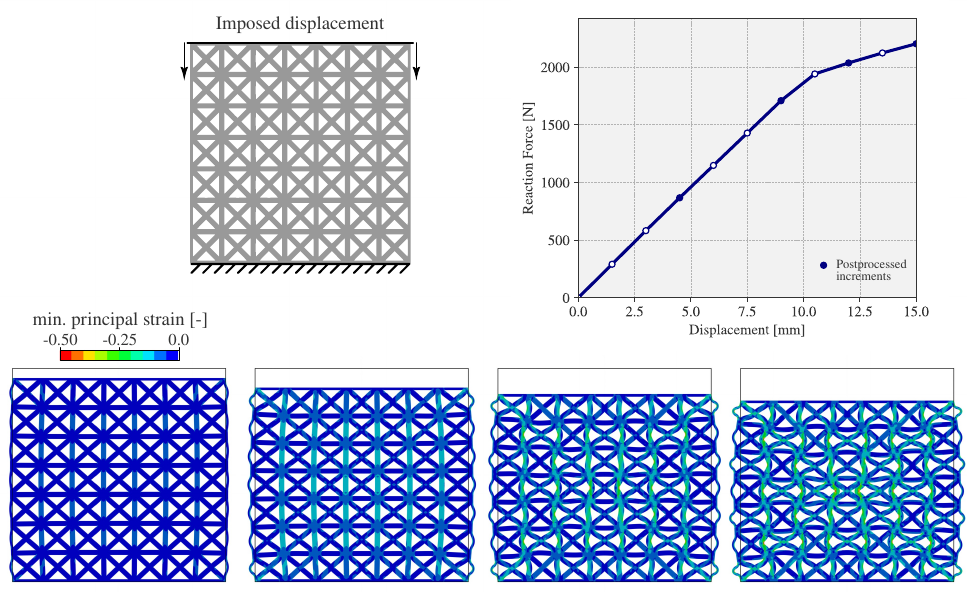}%
	\caption{{Compression of UC1 cells:} an incremental vertical displacement is imposed to the top part leading to $15\%$~compression of the structure. Evolution of the reaction force and the strain field are presented. The lattice structure consists of $N_s=7\times{}7=49$~cells with~$p=3$ and~$h=1/8$.} 
	\label{fignum:8.1}
\end{figure}

\subsubsection{2D straight beam under bending}
Let us start with structural analyses of two-dimensional straight beam lattice structures subjected to bending as illustrated in \cref{fignum:1}. The macroscopic geometry consists of a simple rectangle of dimensions~\qty{100}{\mm} by~\qty{50}{\mm}, modeled by a single linear patch. We employ here the cross-cell UC1. A vertical downward load $g_y = -\qty{3}{\N\per\mm}$ is applied along the right edge of the structure (see again \cref{fignum:1}).

To ease the convergence of the Newton method, the load is applied incrementally in five steps. The deformed configuration of the lattice structure and the associate strain for different number of cells are shown in \cref{fignum:1}. The results exhibit a local buckling phenomenon in some unit cells of the lattice, specifically those located at the bottom-left of the structure. This buckling occurs because these cells experience higher compressive stresses compared to the rest. The buckling allows for load redistribution to neighboring cells. The more the number of cells, the less these local bucklings impact the global mechanical behavior of the structure: the load/displacement curve exhibits less distinct softening.


\begin{figure}[pos=p]
    \centering
    \includegraphics*[width=0.9\linewidth]{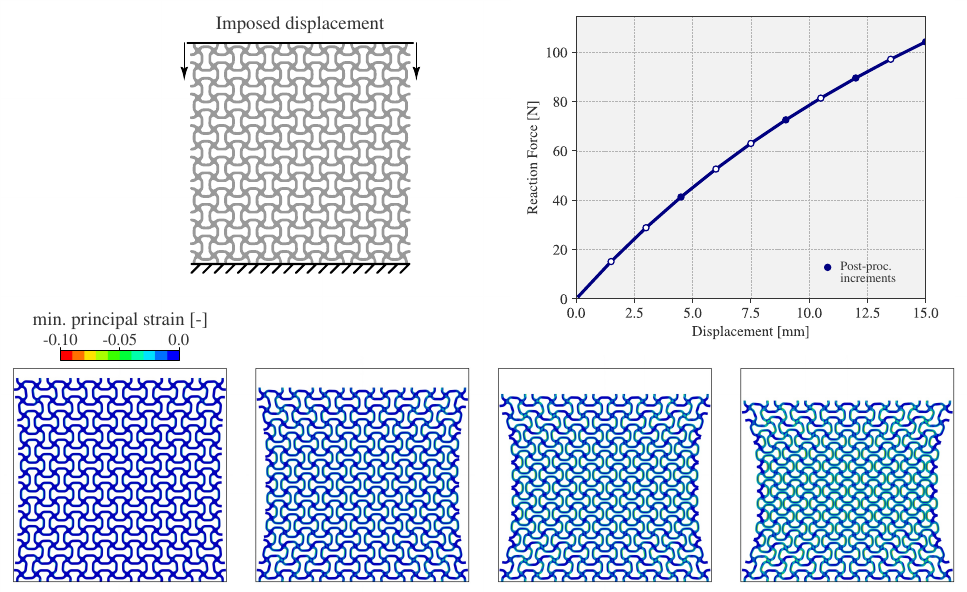}%
	\caption{{Compression of UC2 cells:} an incremental vertical displacement is imposed to the top part leading to $15\%$~compression of the structure. Evolution of the reaction force and the strain field are presented. The lattice structure consists of $N_s=7\times{}7=49$~cells with~$p=3$ and~$h=1/4$.} 
	\label{fignum:8.2}
\end{figure}
\begin{figure}[pos=p]
    \centering
    \includegraphics*[width=0.9\linewidth]{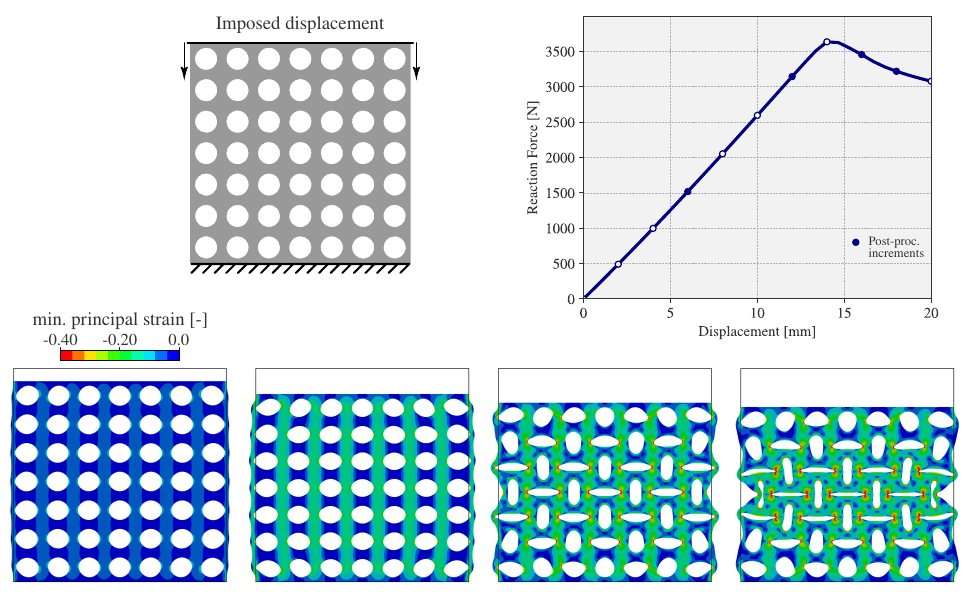}%
	\caption{{Compression of UC3 cells:}  an incremental vertical displacement is imposed to the top part leading to $20\%$~compression of the structure. Evolution of the reaction force and the strain field are presented. The lattice structure consists of $N_s=7\times{}7=49$~cells with~$p=3$ and~$h=1/8$.} 
	\label{fignum:8.3}
\end{figure}
\begin{figure}[pos=t]
    \centering
    \includegraphics[width=\linewidth]{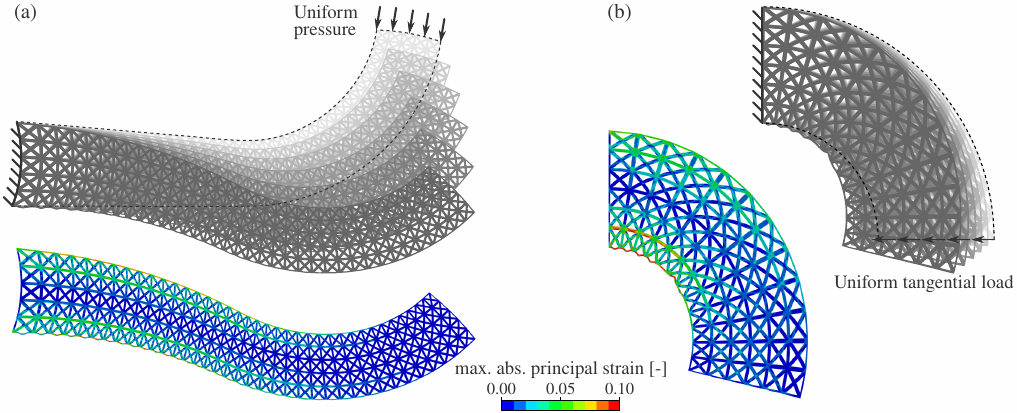}%
	\caption{Additional 2D test cases with UC1 unit cell: (a) brake pedal under operating conditions, and (b) curved beam under bending.} 
	\label{fignum:9}
\end{figure}

\subsubsection{Compression of 2D lattice structures}
As a second test case, we examine the classical compression behavior of lattice structures. The model involves a macroscopic geometry made of a square domain of dimensions $\qty{100}{\mm}\times{}\qty{100}{\mm}$, modeled by a single linear patch. A vertical downward displacement is imposed on the top face of the structure, while horizontal displacement is constrained. The maximal imposed displacement is defined to prevent contact phenomena as contacts are not modeled in this work. We investigate three different unit-cell types (UC1, UC2, and UC3) using degree $p=3$ and mesh refinement $h=1/8$. The results are presented in \cref{fignum:8.1,fignum:8.2,fignum:8.3}, where the strain fields and the load/displacement response function are depicted.


For the UC1 and UC3 unit cells, the response is initially linear up to a compression of approximately $11\%$ strain for UC1 and $14\%$ strain for UC3. Beyond these thresholds, the stress-strain relation becomes (highly) nonlinear, corresponding to the onset of local buckling occurring within the structure. This phenomenon can be clearly observed in the corresponding deformed configurations shown in \cref{fignum:8.1,fignum:8.3}. For the case of cells with circular void (UC3, \cref{fignum:8.3}), we obtain classical deformation patterns, see for instance~\cite{Bertoldi2008,Guo2024,Sperling2024}. As compression increases and the structure exits the linear regime, a pattern transformation occurs, yielding alternating, mutually orthogonal elliptical shapes. As compression continues, the ellipses become increasingly pronounced: their major axes elongate while their minor axes contract with increasing macroscopic compressive strain.

In contrast, the structure composed of UC2 unit cells does not show any local buckling. The mechanical response is stable over the whole studied compression regime as shown in \cref{fignum:8.2}. Finally, such a design exhibits an auxetic behavior under compression: specifically, vertical compression induces horizontal contraction.

\subsubsection{Conformal lattice structures}

Two final two-dimensional configurations, representing a brake pedal lattice structure under operating conditions and a curved beam under bending, are considered. Such lattice structures are often denoted as ``conformal'' lattice structures since the unit cells are not simply repeated (and clipped) with a macro-domain, but they are deformed to match its curved boundaries. These test cases are taken from~\citet{Hirschler:2023aa} where details on the underlying geometric models are provided. The UC1 unit-cell configuration is investigated using degree $p=3$ and mesh refinement $h=1/8$. The resulting maximum principal strain at equilibrium under the applied load is shown in \cref{fignum:9} for both cases. 


\begin{figure}[pos=p]
    \centering
    \includegraphics[width=0.95\linewidth]{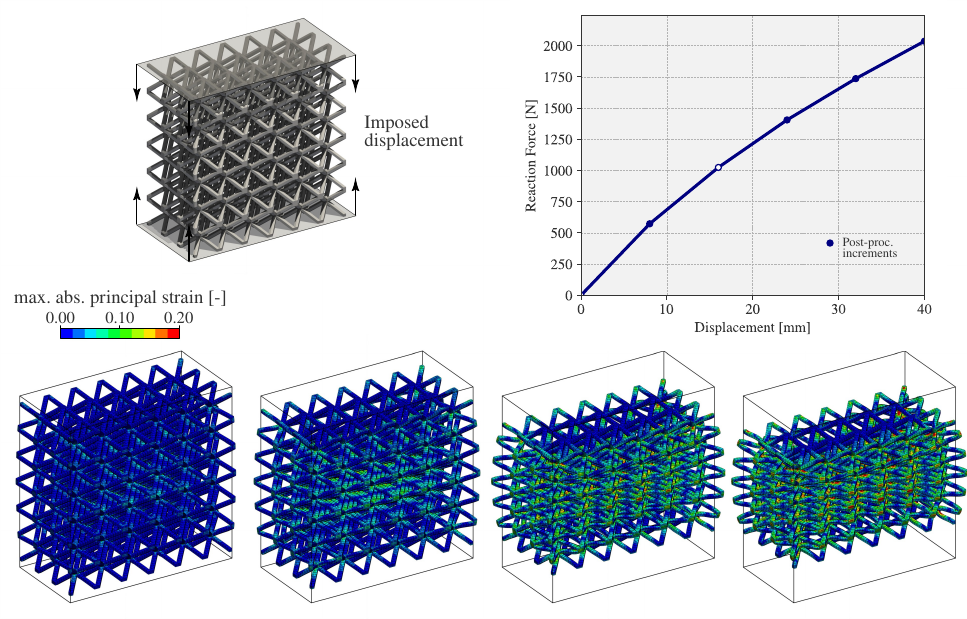}%
	\caption{{3D planar map with UC4 unit cell:} an imposed displacement compresses the structure up to \qty{40}{\percent}. Here the results are given for a total of $N_s=108$~cells with~$p=2$, $h=1/2$.} 
    \label{fignum:7.1}
\end{figure}
\begin{figure}[pos=p]
    \centering
    \includegraphics[width=0.95\linewidth]{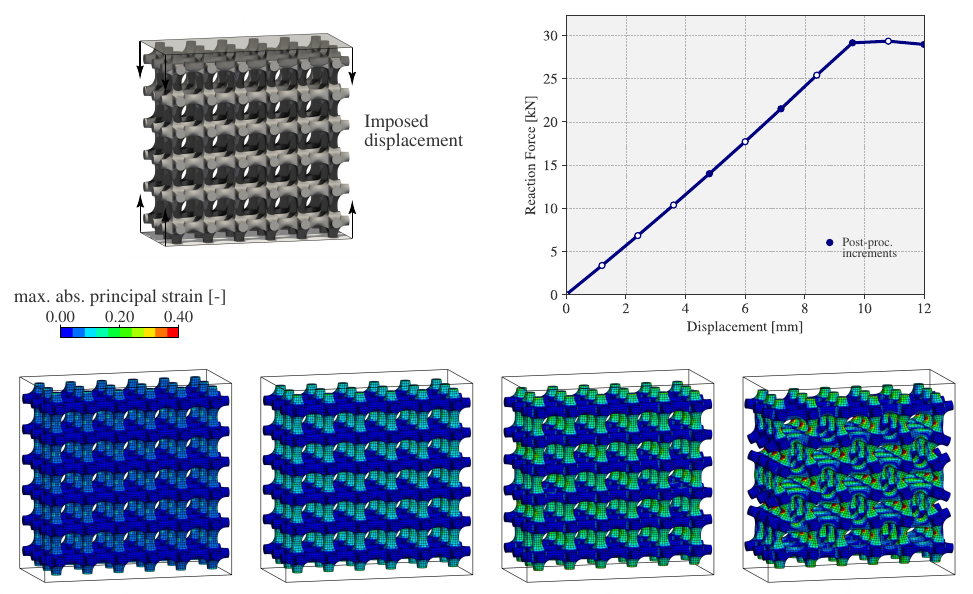}%
	\caption{{3D planar map with UC5 unit cells:} an imposed displacement compresses the structure up to \qty{12}{\percent}. Here the results are given for a total of $N_s=108$~cells with~$p=3$, $h=1/4$.} \label{fignum:7.2}
\end{figure}

\subsubsection{Compression of 3D lattice structures}
Lastly, we consider a three-dimensional configuration involving lattice structures subjected to compression. The macroscopic geometry consists of a hexahedral domain of size~\qtyproduct{100 x 50 x 100}{\mm}, modeled by a single linear patch. Two different unit-cell types (UC4) and (UC5) are investigated, see \cref{fignum:7.1} and \cref{fignum:7.2}. The compression is done via imposed displacement, up to the self-contact starting point. For this test case, the BCC cell (UC4) is stable over the whole deformation regime whereas the cross-cell (UC5) buckles abruptly right before \qty{10}{\percent} compression.

 

\begin{table}[pos=t]
\caption{Computational time and memory footprint for lattice structures with unit cell UC1 ($p=3$ and $h=1/4$) under bending (similar to test case depicted in \cref{fignum:1}) with loads applied over 4 load increments. The RB tolerance is set to $\varepsilon= 3\times10^{-4}$.}
\label{tab:1}

\begin{tabular}{@{}rrrrrr@{}}
\toprule
\multicolumn{6}{l}{\textbf{Standard method}} \\ 
\multicolumn{1}{l}{\#{}Cells} & 
\multicolumn{1}{l}{\#{}DOFs} &
  \multirow{2}{*}{\begin{tabular}[c]{@{}l@{}}Total \\ \#{}Newton iter.\end{tabular}} &
  \multicolumn{2}{l}{Execution time (s)} &
  \multirow{2}{*}{\begin{tabular}[c]{@{}l@{}}Memory\\ (GB)\end{tabular}} \\ \cmidrule(lr){4-5}
 &
   &
   &
  \begin{tabular}[c]{@{}l@{}}Assembly\end{tabular} &
  \begin{tabular}[c]{@{}l@{}}Fact./Solve\end{tabular} &
   \\ \midrule
\num{72} (\numproduct{12 x 06})   & \num{89858}  & \num{20} & \num{108}  & \num{180}  & \num{3.1}  \\ 
\num{392} (\numproduct{28 x 14}) & \num{487426}  & \num{16} & \num{488}  & \num{1940} & \num{6.4}  \\ 
\num{648} (\numproduct{36 x 18}) & \num{805250}  & \num{17} & \num{656}  & \num{3510} & \num{9.6}  \\ 
\num{968} (\numproduct{44 x 22}) & \num{1202434} & \num{19} & \num{1022} & \num{6394} & \num{13.6} \\ 
\midrule
\multicolumn{6}{l}{\textbf{RB method}} \\ 
\multicolumn{1}{l}{\#{}Cells} & 
\multicolumn{1}{l}{\#{}DOFs} &
  \multirow{2}{*}{\begin{tabular}[c]{@{}l@{}}Total \\ \#{}Newton iter.\end{tabular}} &
  \multicolumn{2}{l}{Execution time (s)} &
  \multirow{2}{*}{\begin{tabular}[c]{@{}l@{}}Memory\\ (GB)\end{tabular}} \\ \cmidrule(lr){4-5}
 &
   &
   &
  \multicolumn{1}{r}{Init} &
  \begin{tabular}[c]{@{}l@{}}Solve\end{tabular} &
   \\ \midrule
\num{72} (\numproduct{12 x 06})   & \num{89858}  & \num{20} & \num{66}  & \num{21}  & \num{1.4}  \\ 
\num{392} (\numproduct{28 x 14}) & \num{487426}  & \num{20} & \num{272}  & \num{113} & \num{2.3}  \\ 
\num{648} (\numproduct{36 x 18}) & \num{805250}  & \num{18} & \num{364}  & \num{132} & \num{3.3}  \\ 
\num{968} (\numproduct{44 x 22}) & \num{1202434} & \num{19} & \num{509} & \num{238} & \num{4.2} \\ 
\bottomrule
\end{tabular}
\end{table}

\subsection{Convergence and performance analysis}
\label{sec:num:2}
In this section, we focus on the numerical performance of the proposed method. These include the convergence behavior of the Newton method, the efficiency of the fast-assembly strategy, and the performance of the domain-decomposition solver used at each Newton iteration. In the following, we refer to our approach as the RB method.

To evaluate the performance of the RB method, we compare it against a standard reference method. In this standard approach, the tangent operator is constructed by assembling all the local operators, each computed using full order Gauss quadrature ($p+1$ points per direction per element). The resulting linear system is then solved using a direct LU factorization. \footnote{Note that since the global tangent matrix is symmetric positive definite, a Cholesky factorization could alternatively be used to reduce the computational cost by approximately a factor of two.}. 
The RB method differs from the standard approach in two main aspects. First, in the RB method, the tangent operator is approximated as a linear combination of a small number of local operators, corresponding to principal cells that belong to the RB. Second, the linear system is solved using an inexact FETI-DP method constructed on the RB. We do not apply the RB strategy to the right-hand side, as its computation is inexpensive in our case studies, \corclem{and we require an accurate representation of the residual to ensure that the Newton method has converged to the true solution}. 

A key aspect of the RB method is the construction of an accurate RB at a low computational cost. To achieve this, we employ a greedy algorithm applied to evaluations of the tangent operator at selected points associated with a reduced quadrature rule; see \cref{sec:constr_rb} for the details. 
The key is to employ a quadrature rule that is sufficiently accurate to capture most of the geometry of the macro-mapping and the deformations within the cells, yet inexpensive enough to ensure that the identification process remains significantly cheaper than computing all exact tangent operators. We have observed numerically that an inexpensive quadrature rule is sufficient to accurately identify the principal cells. In fact, using only $\tilde{n}_g = 2^d$ Gauss quadrature points per element proved adequate across all test cases considered in our numerical investigations.

\begin{table}[pos=t]
\caption{Computational time and memory footprint for lattice structures with unit cell UC2 ($p=3$ and $h=1/4$) under compression (similar to test case depicted in \cref{fignum:8.2}), subjected to an imposed displacement leading to \qty{10}{\percent} compression applied over 4 load increments. The RB tolerance is set to $\varepsilon= 3\times10^{-4}$.}
\label{tab:2}
\begin{tabular}{@{}rrrrrr@{}}
\toprule
\multicolumn{6}{l}{\textbf{Standard method}} \\ 
\multicolumn{1}{l}{\#{}Cells} & 
\multicolumn{1}{l}{\#{}DOFs} &
    \multirow{2}{*}{\begin{tabular}[c]{@{}l@{}}Total \\ \#{}Newton iter.\end{tabular}} &
    \multicolumn{2}{l}{Execution time (s)} &
    \multirow{2}{*}{\begin{tabular}[c]{@{}l@{}}Memory\\ (GB)\end{tabular}} \\ \cmidrule(lr){4-5}
 &
     &
     &
    \begin{tabular}[c]{@{}l@{}}Assembly\end{tabular} &
    \begin{tabular}[c]{@{}l@{}}Fact./Solve\end{tabular} &
     \\ \midrule
\num{64} (\numproduct{8 x 8})      & \num{16962}  & \num{22} & \num{30}  & \num{34}   & \num{0.57} \\
\num{256} (\numproduct{16 x 16})   & \num{67458}  & \num{16} & \num{93}  & \num{262}  & \num{1.4}  \\
\num{576} (\numproduct{24 x 24})   & \num{151490} & \num{16} & \num{19}  & \num{778}  & \num{3.4}  \\
\num{1024} (\numproduct{32 x 32})  & \num{269058} & \num{16} & \num{330} & \num{1853} & \num{6.3}  \\
\midrule
\multicolumn{6}{l}{\textbf{RB method}} \\ 
\multicolumn{1}{l}{\#{}Cells} & 
\multicolumn{1}{l}{\#{}DOFs} &
    \multirow{2}{*}{\begin{tabular}[c]{@{}l@{}}Total \\ \#{}Newton iter.\end{tabular}} &
    \multicolumn{2}{l}{Execution time (s)} &
    \multirow{2}{*}{\begin{tabular}[c]{@{}l@{}}Memory\\ (GB)\end{tabular}} \\ \cmidrule(lr){4-5}
 &
     &
     &
    \multicolumn{1}{r}{Init} &
    \begin{tabular}[c]{@{}l@{}}Solve\end{tabular} &
     \\ \midrule
\num{64} (\numproduct{8 x 8})      & \num{16962}  & \num{15} & \num{16}  & \num{1}    & \num{0.3} \\
\num{256} (\numproduct{16 x 16})   & \num{67458}  & \num{16} & \num{49}  & \num{7}    & \num{0.8} \\
\num{576} (\numproduct{24 x 24})   & \num{151490} & \num{16} & \num{129} & \num{25}   & \num{1.3} \\
\num{1024} (\numproduct{32 x 32})  & \num{269058} & \num{16} & \num{260} & \num{76}   & \num{2.2} \\
\bottomrule
\end{tabular}
\end{table}

\begin{table}
\caption{Computational time and memory footprint for lattice structures with BCC-cell UC4 ($p=1$ and $h=1/2$) under compression (similar to test case depicted in \cref{fignum:7.1}), subjected to an imposed displacement applied over 4 load increments leading to \qty{10}{\percent} compression. The RB tolerance is set to $\varepsilon=5\times 10 ^{-3}$. $\dagger$ means out of memory.}
\label{tab:3}
\begin{tabular}{@{}rrrrrr@{}}
\toprule
\multicolumn{6}{l}{\textbf{Standard method}} \\ 
\multicolumn{1}{l}{\#{}Cells} & 
\multicolumn{1}{l}{\#{}DOFs} &
    \multirow{2}{*}{\begin{tabular}[c]{@{}l@{}}Total \\ \#{}Newton iter.\end{tabular}} &
    \multicolumn{2}{l}{Execution time (s)} &
    \multirow{2}{*}{\begin{tabular}[c]{@{}l@{}}Memory\\ (GB)\end{tabular}} \\ \cmidrule(lr){4-5}
 &
     &
     &
    \begin{tabular}[c]{@{}l@{}}Assembly\end{tabular} &
    \begin{tabular}[c]{@{}l@{}}Fact./Solve\end{tabular} &
     \\ \midrule
\num{32} (\numproduct{4 x 4 x 2})    & \num{21501}   & \num{13}          & \num{9.6}   & \num{25}    & \num{0.6}  \\
\num{256} (\numproduct{8 x 8 x 4})   & \num{167991}  & \num{13}          & \num{74}    & \num{1878}  & \num{5.4}  \\
\num{864} (\numproduct{12 x 12 x 6})  & \num{526737}  & \num{12}          & \num{166}   & \num{12704} & \num{11.3} \\
\num{2048} (\numproduct{16 x 16 x 8}) & \num{1329003} & {$\dagger$} & {$\dagger$} & {$\dagger$} & $\dagger$ \\
\midrule
\multicolumn{6}{l}{\textbf{RB method}} \\ 
\multicolumn{1}{l}{\#{}Cells} & 
\multicolumn{1}{l}{\#{}DOFs} &
    \multirow{2}{*}{\begin{tabular}[c]{@{}l@{}}Total \\ \#{}Newton iter.\end{tabular}} &
    \multicolumn{2}{l}{Execution time (s)} &
    \multirow{2}{*}{\begin{tabular}[c]{@{}l@{}}Memory\\ (GB)\end{tabular}} \\ \cmidrule(lr){4-5}
 &
     &
     &
    \multicolumn{1}{r}{Init} &
    \begin{tabular}[c]{@{}l@{}}Solve\end{tabular} &
     \\ \midrule
\num{32} (\numproduct{4 x 4 x 2})    & \num{21501}   & \num{17} & \num{38}    & \num{4.1}  & \num{0.9}  \\
\num{256} (\numproduct{8 x 8 x 4})   & \num{167991}  & \num{20} & \num{385}   & \num{89}   & \num{2.7}  \\
\num{864} (\numproduct{12 x 12 x 6})  & \num{526737}  & \num{20} & \num{1952}  & \num{447}  & \num{5.4}  \\
\num{2048} (\numproduct{16 x 16 x 8}) & \num{1329003} & \num{23} & \num{10289} & \num{2265} & \num{10.5} \\
\bottomrule
\end{tabular}
\end{table}

In \cref{tab:1,tab:2,tab:3}, we compare the standard and RB methods across three representative test cases: a 2D lattice structure with UC1 unit cells under bending, a 2D lattice structure with UC2 unit cells under compression, and a 3D lattice structure with UC4 unit cells under compression.
In all simulations, the number of load increments is set to four, and the RB tolerance is fixed to $\varepsilon = 3\times10^{-4}$ (2D cases) or $\varepsilon = 5\times10^{-3}$ (3D case). For each test case, we consider configurations with an increasing number of unit cells, while keeping the mesh size $h$ and the spline degree $p$ fixed.
For both methods, the tables report several key quantities: the total number of Newton iterations (i.e., the sum of Newton iterations over the four load increments), the execution time, and the memory footprint. The execution time is decomposed into two main contributions. For the standard method, the first corresponds to the local computations and the assembly of the global tangent matrix, and the second to the LU factorization and solution of the resulting linear system.
For the RB method, the first contribution corresponds to the initialization phase, which includes both the construction of the RB and the setup of the FETI-DP operators for all principal cells. Note that these two steps are performed at each Newton iteration. The second contribution corresponds to the solution of the linear system using the inexact FETI-DP solver.

In \cref{tab:1,tab:2}, we observe that the number of Newton iterations is comparable between the standard and RB methods across all configurations. This indicates that the reduced basis is sufficiently rich and that the tangent operator approximation is accurate enough for the quasi-Newton method to converge. Both memory usage and computational time are reduced relative to the standard method, with memory requirements decreasing by a factor of three and computational time by nearly an order of magnitude for the largest configuration. 

In \cref{tab:3}, the RB method requires slightly more Newton iterations than the standard method. This is due to the RB tolerance being set to $\varepsilon=5\times 10^{-3}$, larger than in the previous test cases. While this choice reduces memory usage, it increases computational time because more Newton iterations are needed. We note that this choice may not represent the optimal trade-off between memory and computational cost. Nevertheless, even with this setting, the RB method is approximately five times faster than the standard method. Additionally, we successfully performed a simulation with the RB method in a configuration that exceeded the memory capacity of the standard method.

\section{Conclusions}
\label{sec:conclu}

In this work, we developed a method for the efficient fine-scale simulation of nonlinear hyperelastic lattice structures, also referred to as architected or metamaterials. Unlike existing approaches, which typically avoid a full volumetric fine-scale analysis, either through multiscale homogenization schemes or beam and shell network models, the strategy proposed here starts with a complete three-dimensional discretization of the entire lattice architecture (thus including all struts and/or walls), and then designs a dedicated solver that drastically reduces both memory storage and computational cost. Building upon previous contributions that focused on the linear elastic regime~\cite{Hirschler2022,Hirschler:2023aa}, the key idea was to exploit the intrinsic geometrical nature of such materials, namely their natural non-overlapping domain decomposition into a collection of similar cells. 

To effectively extend the approach to the nonlinear setting, we designed a ROM approach that enables the efficient identification of so-called principal cells throughout the Newton iterations. In practice, this consisted in extracting, at each Newton iteration, a limited set of principal cells; that is, for a given deformation level within the structure, a small number of representative cells whose mechanical behaviors capture, to a large extent, that of all the others. 
As a first contribution, this resulted in fast and memory-efficient operator assembly strategies, in which only the operators of the principal cells were fully assembled. 
This ROM procedure can thus be interpreted as an EIM, where the (predefined) ``magic points'' correspond to the reduced integration points. Secondly, the benefit of the fast assembly approach was extended to the solution stage of the nonlinear problem. The constructed ROM was used to precondition a domain-decomposition solver, where each cell is regarded as a subdomain. More precisely, the principal tangent matrices were employed within the inexact FETI-DP algorithm proposed in~\cite{Hirschler:2023aa}, which significantly reduces the number of required local factorizations. Consequently, we ended up with a quasi matrix-free algorithm that drastically lowers the computational time for solving the tangent linear systems at each Newton iteration. 

Several two- and three-dimensional numerical experiments, covering various cell patterns and macro-mappings, were performed to assess the performance of the proposed solver. The simulations were run up to deformation levels of about 30–40\%, corresponding to the state of the structure just before the onset of self-contact. At this stage, only serial computations were conducted on a standard laptop. Comparisons with a conventional direct solver were also carried out to highlight the advantages of the proposed approach. Overall, the computational time was reduced from several hours to a few tens of minutes, while the memory usage decreased by a factor of about three. This notably made it possible to compute problems involving thousands of cells (i.e., millions of degrees of freedom) within a few minutes on an off-the-shelf laptop. The encouraging results obtained in this study motivate further developments along this line for nonlinear simulations, aiming to exploit the self-similarity of lattice cells not at the modeling stage but at the solver level. Future work could focus on extending the present framework to elastoplasticity and self-contact, as well as on coupling the proposed approach with dedicated high-performance computing algorithms, as in~\cite{GUILLET2025114136}, in order to enable the simulation of lattice structures composed of several tens of thousands of cells. Regarding this last point, exploring parallel-in-Newton domain decomposition strategies (see, e.g.,~\cite{cai2002,dolean2016nonlinear,klawonn17}) may constitute a promising direction. 

\paragraph*{Acknowldegments:}

This research was supported by the French ``Agence Nationale de la Recherche'' under grant ANR-22-CE46-0007 (AVATAR). ANR is gratefully acknowledged. As part of the ``France 2030'' initiative, this work has benefited from a national grant managed by ``Agence Nationale de la Recherche'' attributed to the Exa-MA project of the NumPEx PEPR program, under the reference ANR-22-EXNU-0002. Pablo Antolin acknowledges the financial support of the Swiss National Science Foundation through the project FLAS$_h$ with no.\ 200021\_214987.
For the purpose of Open Access, a \href{https://creativecommons.org/licenses/by/4.0/}{[CC-BY public copyright license]} has been applied by the authors to the present document and will be applied to all subsequent versions up to the Author Accepted Manuscript arising from this submission.

\paragraph*{Declaration of generative AI and AI-assisted technologies in the writing process:}

During the preparation of this work, the authors used ChatGPT in order to improve language and readability.  After using this tool/service,  the authors reviewed and edited the content as needed and take full responsibility for the content of the publication.

\bibliography{bib/bib}
\bibliographystyle{cas-model2-names}
\end{document}